\documentclass{newtonpreprint1}
\usepackage[dvips]{graphicx}
\usepackage{graphics} 
\usepackage{epsf} 
\usepackage{amsmath, amsfonts, amssymb, graphicx, color, subfigure, wrapfig, wasysym}

 \theoremstyle{definition}
 
 \theoremstyle{remark}

\begin{document}

\title[Fixed equilibria]{Equilibrium singularity distributions in the plane}

\author[P.K. Newton \& V. Ostrovskyi]{Paul K. Newton and Vitalii Ostrovskyi}

\affiliation{Department of Aerospace \& Mechanical Engineering and Department of Mathematics\\
University of Southern California, Los Angeles, CA 90089-1191}

\label{firstpage}

\maketitle

\begin{abstract}{Singularity dynamics; Fixed equilibria; Singular values; Point vortex equilibria; Shannon entropy}
We characterize  all {\it fixed} equilibrium point singularity distributions in the plane of logarithmic type, allowing for real, imaginary,
or complex singularity strengths $\vec{\Gamma}$. The dynamical system  follows from the assumption that each
of the $N$ singularities  moves according to the flowfield generated by all the others at that point.
For  strength vector $\vec{\Gamma} \in {\Bbb R}^N$, the dynamical system is the classical
point vortex system obtained from a singular discrete representation of the vorticity field from incompressible fluid flow. When $\vec{\Gamma} \in Im$,  it corresponds to a system of sources and sinks,
whereas when $\vec{\Gamma} \in {\Bbb C}^N$ the system consists of spiral sources and sinks discussed
in Kochin et. al.  (1964). We formulate the equilibrium problem as one in linear algebra,
$A \vec{\Gamma} = 0$, $A \in {\Bbb C}^{N \times N}$, $\vec{\Gamma} \in {\Bbb C}^N$, where
$A$ is  a $N \times N$ complex skew-symmetric configuration matrix which encodes the geometry of
the system of interacting singularities. For an equilibrium to exist,  $A$ must have a
kernel.
$\vec{\Gamma}$ must  then be an element of the nullspace
of $A$.
We prove that when
$N$ is odd, $A$ always has a kernel,  hence there is a choice of $\vec{\Gamma}$ for which
the system is a fixed equilibrium. When $N$ is even, there may or may not be a non-trivial
nullspace of $A$, depending on the relative position of the points in the plane.
We  describe a  method for classifying the equilibria in terms of the distribution
of the non-zero eigenvalues (singular values) of $A$, or equivalently, the non-zero eigenvalues of the
associated covariance matrix $A^{\dagger} A$, from which one can calculate the Shannon entropy of the
configuration.
\end{abstract}

\section{Introduction}

Consider the vector field at $z = 0$ governed by the complex dynamical system:
\begin{eqnarray}
\dot{z}^* = \frac{\Gamma}{2 \pi i} \frac{1}{z}, \quad z(t) \in {\Bbb C}, \quad   \Gamma \in {\Bbb C}, \quad
t \in {\Bbb R} > 0,\label{z}
\end{eqnarray}
where $z^*$ denotes the complex conjugate of $z(t)$.
Letting $z(t) = r(t) \exp (i \theta(t))$, $\Gamma = \Gamma_r + i \Gamma_i$, gives:
\begin{eqnarray}
\dot{r} &=& \frac{\Gamma_i}{2 \pi r},\\
\dot{\theta} &=& \frac{\Gamma_r}{2 \pi r^2},
\end{eqnarray}
from which it is easy to see that:
\begin{eqnarray}
r(t) &=& \sqrt{\left(\frac{\Gamma_i}{2\pi}\right) t + r^2 (0)},
\end{eqnarray}

\begin{equation}
\theta(t) = \left\{ \begin{array}{ll}
\left(\frac{\Gamma_r}{\Gamma_i} \right) \ln \left( \left(\frac{\Gamma_r}{\Gamma_i}\right) t + r^2 (0)\right) & \mbox{if $ \Gamma_i \ne 0$} \\\\
\frac{\Gamma_r t}{2 \pi r^2 (0)} + \theta (0) & \mbox{if $ \Gamma_i = 0$}.
\end{array}
\right.
\end{equation}

\begin{figure}[ht]
\begin{tabular}{cccc}
\includegraphics[scale=0.25,angle=0]{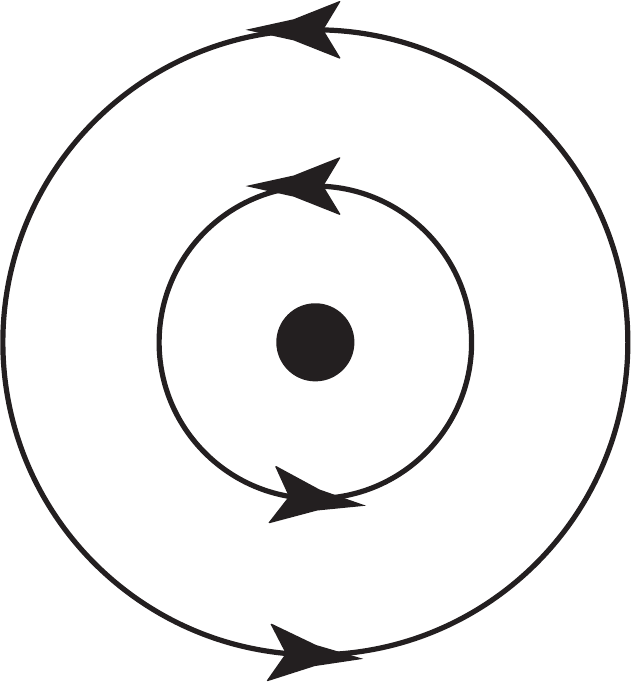} &
\includegraphics[scale=0.25,angle=0]{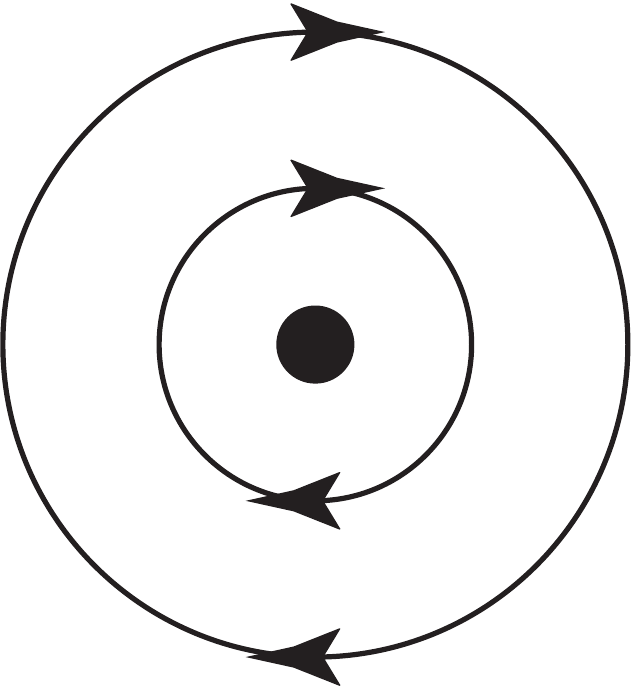}&
\includegraphics[scale=0.25,angle=0]{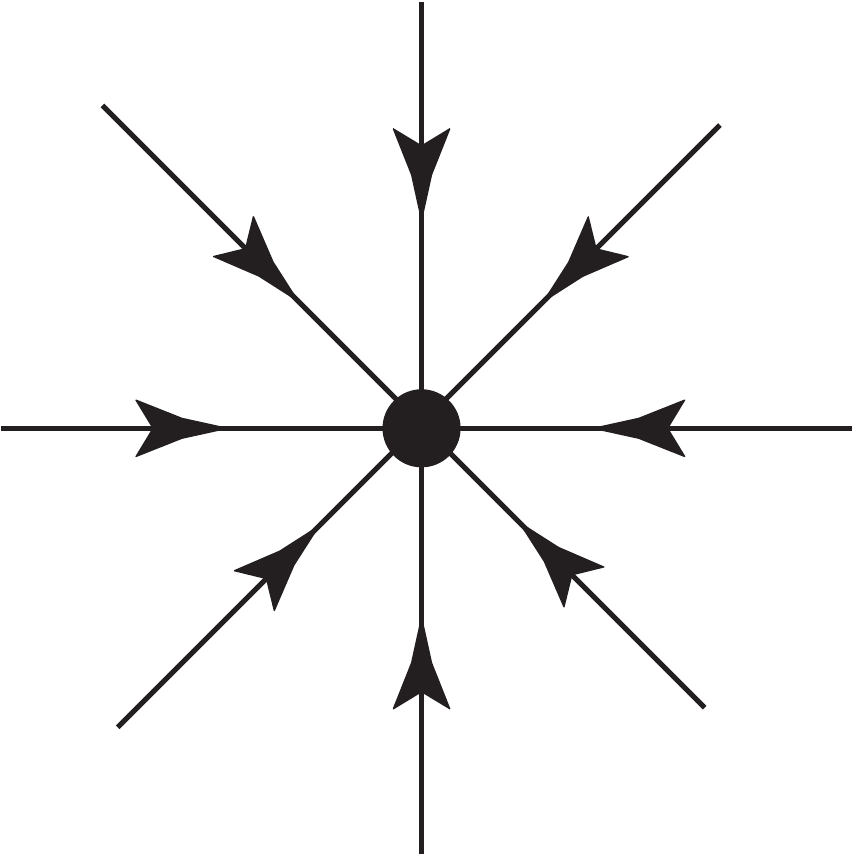} &
\includegraphics[scale=0.25,angle=0]{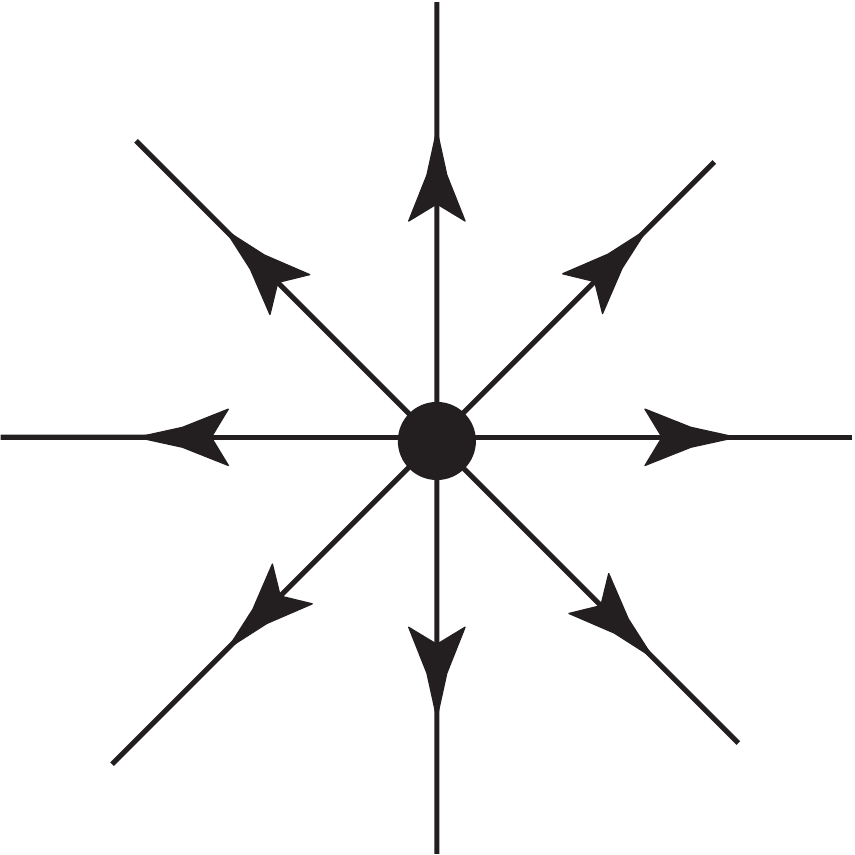}\\
{\footnotesize (a) $\Gamma_r  > 0, \Gamma_i = 0$} & {\footnotesize (b) $\Gamma_r  < 0, \Gamma_i = 0$ } & {\footnotesize (c) $\Gamma_r = 0, \Gamma_i   < 0$} & {\footnotesize (d) $\Gamma_r = 0, \Gamma_i > 0$}\\
\vspace{0.01cm}\\
\includegraphics[scale=0.25,angle=0]{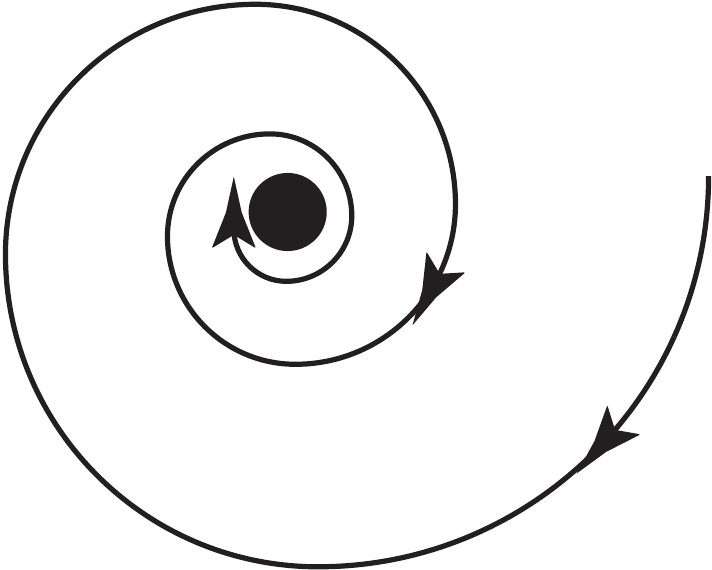} &
\includegraphics[scale=0.25,angle=0]{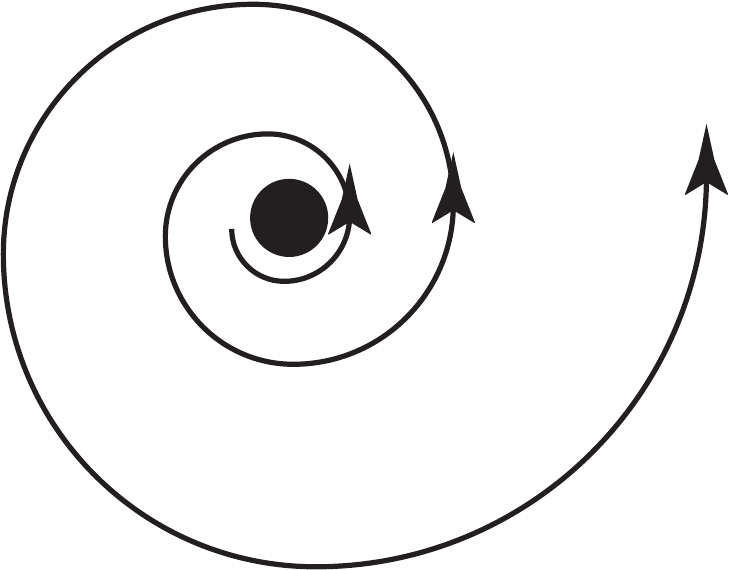}&
\includegraphics[scale=0.25,angle=0]{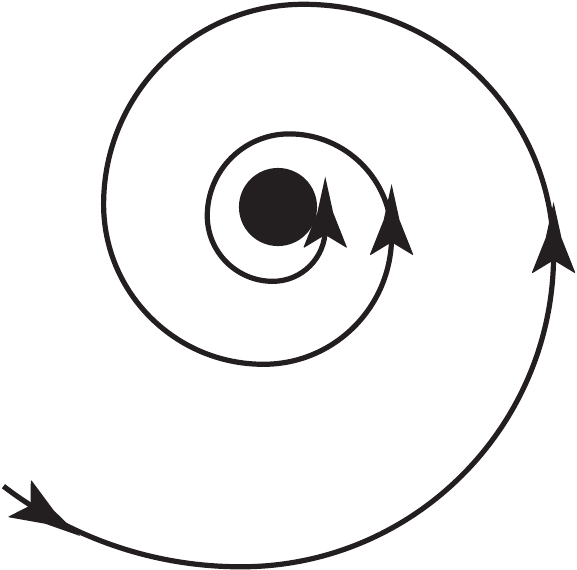} &
\includegraphics[scale=0.25,angle=0]{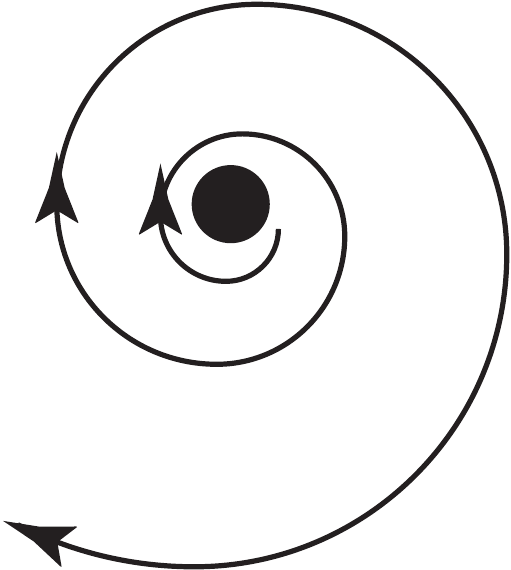}\\
{\footnotesize (e) $\Gamma_r  < 0$, $ \Gamma_i  < 0$} & {\footnotesize (f)  $ \Gamma_r  > 0$, $\Gamma_i  > 0$} & {\footnotesize (g)  $\Gamma_r  > 0$, $ \Gamma_i  < 0$ } & {\footnotesize (h)  $\Gamma_r  < 0$, $\Gamma_i  >  0$}\\
\vspace{0.01cm}\\
\end{tabular}
\caption{\footnotesize
All possible flowfields at the singular point $z = 0$ associated with the dynamical system (\ref{z}).}
\label{fig1}
\end{figure}

\noindent
When $\Gamma_r \ne 0$, $\Gamma_i = 0$,  the field is that of a classical point-vortex (figure \ref{fig1}(a),(b));
when $\Gamma_r = 0$, $\Gamma_i \ne 0$ it is a source ($\Gamma_i > 0$) or sink ($\Gamma_i < 0$)
(figure \ref{fig1}(c),(d)), while when
$\Gamma_r \ne 0$, $\Gamma_i \ne 0$, it is a spiral-source or sink
((figure \ref{fig1}(e)-(h)).

A  collection of $N$ of these point singularities, each located at $z = z_{\beta}(t)$, $\beta = 1, ..., N$,
by linear superposition, produces the field:
\begin{equation}
{\dot{z}^{*}} = \frac{1}{2 \pi i} \sum_{\beta = 1}^N  \frac{\Gamma_{\beta}}{z - z_{\beta}};
\quad z(t) \equiv x(t) + iy(t) \in {\Bbb C},\quad  \Gamma_{\beta} \in {\Bbb C}.\label{eqn1}
\end{equation}
Then, if we  advect
each  by the velocity field generated by all the others\footnote{One might characterize this dynamical assumption by saying that each singularity `goes with the flow'.}, we arrive at the complex dynamical system:
\begin{equation}
{\dot{z}^{*}_{\alpha}} = \frac{1}{2 \pi i} \sum_{\beta = 1}^N{}' \frac{\Gamma_{\beta}}{z_{\alpha} - z_{\beta}};
\quad z_{\alpha}(t) \equiv x_{\alpha}(t) + iy_{\alpha}(t) \in {\Bbb C},\quad  \Gamma_{\beta} \in {\Bbb C},\label{eqn1a}
\end{equation}
where $'$ indicates  that $\beta \ne \alpha$.
In this paper we  characterize  all {\it fixed} equilibria of (\ref{eqn1a}), namely solutions
for which ${\dot{z}^{*}}_{\alpha}(t) = 0$. For this, we have the $N$ coupled equations:
\begin{equation}
\sum_{\beta = 1}^N{}' \frac{\Gamma_{\beta}}{z_{\alpha} - z_{\beta}} = 0, \quad
(\alpha = 1, ... N),\label{eqn2}
\end{equation}
where we are interested in positions $z_{\alpha}$ and strengths $\Gamma_{\alpha}$ for which this nonlinear algebraic
system is satisfied.
Since Eqn (\ref{eqn2})   is {\it linear} in the $\Gamma$'s,  it can more productively  be written in matrix form
\begin{equation}
A \vec{\Gamma} = 0\label{eqn3}
\end{equation}
where $A \in {\Bbb C}^{N \times N}$ is evidently a skew-symmetric matrix $A = -A^T$, with entries
$[ a_{\alpha \alpha } ] = 0$, $[ a_{\alpha \beta} ] = \frac{1}{z_\alpha - z_\beta} = - [ a_{\beta \alpha} ] $.
We call $A$ the {\it configuration matrix} associated with the interacting particle system (\ref{eqn1a}).
The collection of points $\{ z_1 (0), z_2 (0), ... , z_N (0) \}$
 in the complex plane  is called the {\it configuration}.
%
From (\ref{eqn3}), we can conclude that the points $z_{\alpha}$ are in a fixed equilibrium configuration if
$\det (A) = 0$, i.e. there is at least one zero eigenvalue of $A$.
If the corresponding eigenvector is real, the configuration is made up of point-vortices.
If it  is imaginary, it is made up of sources and sinks.
If it is complex, it is made up of spiral sources and sinks. Notice also that if
$\frac{d{z}^{*}_{\alpha}}{dt} = 0$, then one can prove that
$\frac{d^n{z}^{*}_{\alpha}}{dt^n} = 0$ for any $n$.
It follows  that:

\medskip

\noindent
{\bf Theorem 1.} {\it For a given configuration of $N$ points  $\{z_1 , z_2, ..., z_N \}$ in the complex plane, there exists a  set of singularity strengths $\vec{\Gamma}$
for which the configuration is a fixed equilibrium solution of the dynamical system (\ref{eqn1a})
iff  $A$ has a kernel,
 or equivalently, if there is at least one zero eigenvalue of $A$.  If the nullspace dimension of
 $A$ is one, i.e. there is only one zero eigenvalue, the choice of $\vec{\Gamma}$ is unique (up to a multiplicative constant).
 If the nullspace dimension is greater than one, the choice of $\vec{\Gamma}$ is not unique and can be any linear
 combination of the basis elements of null($A$).}

 The equilibria we consider in this paper all have one-dimensional nullspaces and odd $N$. The more delicate cases
 of equilibria with higher dimensional nullspaces and even $N$ are deferred to a separate study.
We mention here work of Campbell \& Kadtke (1987) and Kadtke \& Campbell (1987) in which a different technique
is described to find stationary solutions to (\ref{eqn1a}).

\section{General properties of the configuration matrix}

Since  $A$ is skew-symmetric, it follows that
\begin{equation}
\det (A) =  \det (- A^T ) = (-1)^N \det (A^T ) = \det (A^T ).\label{eqn4}
\end{equation}
Hence,  for $N$ {\it odd}, we have $- \det (A^T ) = \det (A^T)$, which implies
$\det (A^T ) = 0$.
\medskip

\noindent
{\bf Theorem 2.} {\it When $N$ is odd,
$A$ always has at least one zero eigenvalue, hence
for any configuration there exists  a choice $\vec{\Gamma} \in {\Bbb C}$ for which the
system is a fixed equilibrium.}
\medskip

\noindent
When $N$ is even, there may or may not be a fixed equilibrium, depending on whether or not $A$ has a non-trivial nullspace. In general, we would like to determine a basis set for the nullspace of $A$ for a given configuration, i.e. the set of all strengths for which a given configuration remains fixed.
Other important general properties of skew-symmetric matrices are listed below:
\begin{enumerate}
\item  The eigenvalues always come in pairs $\pm \lambda$. If $N$ is odd, there is one unpaired eigenvalue that is zero.
\item  If $N$ is even, $\det (A) = Pf(A) ^2 \ge 0$, where $Pf$ is the Pfaffian.
\item  Real skew-symmetric matrices have pure imaginary eigenvalues.
\end{enumerate}

Recall that
every matrix can be written as the sum of a Hermitian matrix ($B = B^\dagger$) and a skew-Hermitian
matrix ($C = -C^\dagger$). To see this, notice
\begin{equation}
A \equiv \frac{1}{2}(A + A^{\dagger}) + \frac{1}{2}(A - A^{\dagger}).
\end{equation}
Here, $B \equiv  \frac{1}{2}(A + A^{\dagger}) = B^{\dagger}$ and $C \equiv \frac{1}{2}(A - A^{\dagger}) = -C^{\dagger}$.
A matrix is {\it normal} if $AA^\dagger = A^\dagger A$, otherwise it is {\it non-normal}.
If we calculate $AA^\dagger - A^\dagger A$, where $A = B + C$ as above, then it is easy to see that
\begin{equation}
AA^\dagger - A^\dagger A = 2(CB - BC).
\end{equation}
Therefore, if $B = 0$ or $C = 0$, $A$ is normal.
\medskip

\noindent
{\bf Theorem 3.} {\it All Hermitian or skew-Hermitian matrices are normal.}
\medskip

The generic configuration matrix $A$ arising from (\ref{eqn3}) is, however, non-normal.

\subsection{Spectral decomposition of normal and non-normal matrices}

For normal matrices, the following spectral-decomposition holds:
\medskip

\noindent
{\bf Theorem 4.} {\it  $A$ is a normal  matrix $\Leftrightarrow$ $A$ is unitarily diagonalizable, i.e.
\begin{equation}
A = Q \Lambda Q^\dagger
\label{diag}
\end{equation}
where $Q$ is unitary.}
\medskip

\noindent
Here, the columns of $Q$ are the $N$ linearly independent eigenvectors of $A$ that can be made mutually orthogonal. The matrix $\Lambda$ is a diagonal matrix with  the $N$ eigenvalues down the diagonal. See Golub and Van Loan (1996) for details.

In general, however,   for the system of interacting particles governed by
(\ref{eqn2}), (\ref{eqn3}), $A \in {\Bbb C}^{N \times N}$ will be a non-normal matrix. The most comprehensive decomposition of $A$ in this
case is the singular value decomposition (Golub and Van Loan (1996), Trefethen and Bau (1997)).
 It is a factorization that greatly generalizes the spectral
decomposition of a normal matrix, and it is available for any matrix.

The $N$ singular values, $\sigma^{(i)}$ ($i = 1, \ldots N$), of  $A$, are non-negative real
numbers that satisfy
\begin{align}
A {\bf v}^{(i)} = \sigma^{(i)} {\bf u}^{(i)}; \quad A^\dagger {\bf u}^{(i)} = \sigma^{(i)} {\bf v}^{(i)},\label{singeqns}
\end{align}
where ${\bf u}^{(i)} \in {\Bbb C}^N$ and ${\bf v}^{(i)} \in {\Bbb C}^N$.
The vector ${\bf u}^{(i)}$ is called the left-singular vector associated with $\sigma^{(i)}$, while
${\bf v}^{(i)}$ is the right-singular vector. In terms of these, the matrix $A$ has the factorization
\begin{align}
A = U\Sigma V^\dagger =
\sum_{i=1}^k \sigma^{(i)} {\bf u}^{(i)} {\bf v}^{(i)}{}^T, \quad (k \le N)
\label{svddecomp}
\end{align}
where $U \in {\Bbb C}^{N \times N}$ in unitary, $V \in {\Bbb C}^{N \times N}$ is unitary, and
$\Sigma \in {\Bbb R}^{N \times N}$ is diagonal.
(\ref{svddecomp})  is the non-normal analogue of the spectral decomposition formula (\ref{diag})
where the summation term on the right hand side gives an (optimal) representation of $A$ as a linear combination of
rank-one matrices with weightings governed by the singular values ordered from largest to smallest.
Here, the rank of $A$ is $k$.
The columns of $U$ are the left-singular vectors ${\bf u}^{(i)}$,
while the columns of $V$ are the right-singular vectors ${\bf v}^{(i)}$.
The matrix $\Sigma$ is given by:
\begin{align}
\Sigma = \left( \begin{array} {ccc}\sigma^{(1)} & \cdots  & 0\\ & \ddots & \\ 0 & \cdots & \sigma^{(N)}
  \end{array} \right) \in {\Bbb R}^{N \times N}.\label{singmatrix}
\end{align}
The singular values can be ordered so that
$\sigma^{(1)} \ge \sigma^{(2)} \ge \ldots \ge \sigma^{(N)} \ge 0$
and  one or more may be zero.
As is evident from multiplying the first equation in  (\ref{singeqns}) by $A^\dagger$ and the second
by $A$,
\begin{align}
(A^{\dagger} A - \sigma^{(i)}{}^2 ){\bf v}^{(i)} = 0; \quad (A A^{\dagger} - \sigma^{(i)}{}^2 ){\bf u}^{(i)} = 0,\label{eigs}
\end{align}
the singular values squared are the eigenvalues of the {\it covariance} matrices $A^{\dagger} A$ or $A A^{\dagger}$,
which have the same eigenvalue structure, while the left-singular vectors ${\bf u}^{(i)}$ are the
eigenvectors of $AA^\dagger$, and the right-singular vectors ${\bf v}^{(i)}$ are the eigenvectors of
$A^{\dagger} A$.
From (\ref{singeqns}), we also note that the right singular vectors ${\bf v}^{(i)}$ corresponding to
$\sigma^{(i)} = 0$ form a basis for the nullspace of $A$.
Because of (\ref{eqn3}), we seek configuration matrices with one or more singular values that are zero.

\section{Collinear equilibria}

\begin{figure}[h]
\begin{center}
\includegraphics[scale=0.6,angle=0]{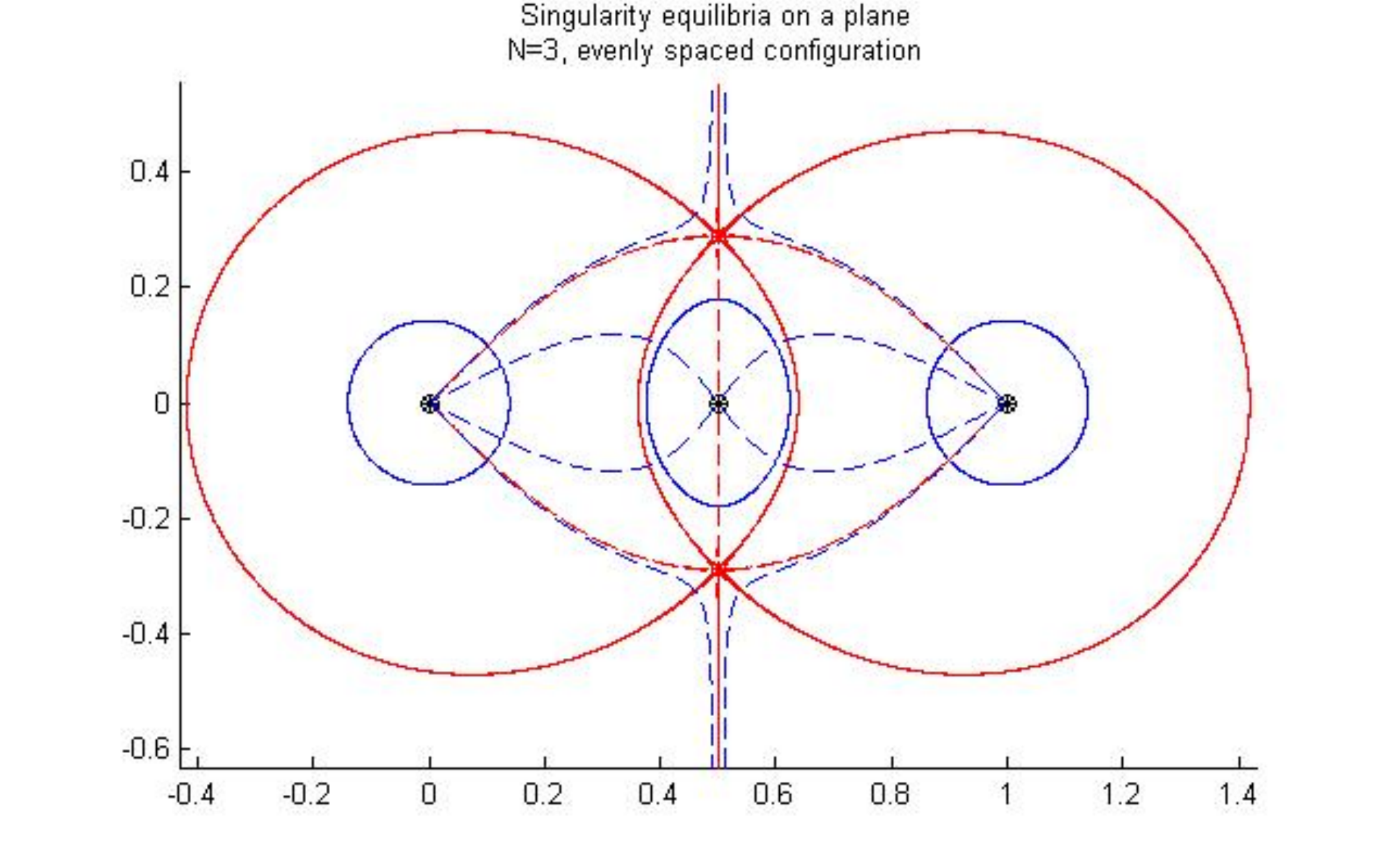}
\end{center}
\caption{\footnotesize  $N = 3$ evenly distributed point vortices  on a line with strengths
$\Gamma_1 = 1, \Gamma_2 = -\frac{1}{2}, \Gamma_3 = 1$, in equilibrium. The far field is that of a point vortex at the
center-of-vorticity of the system. Solid streamline pattern is for point vortices, dashed streamline pattern is for source/sink system. The patterns are orthogonal.}
\label{fig2}
\end{figure}

For the special case in which all the particles lie on a straight line, there is no loss in assuming
$z_{\alpha} = x_{\alpha} \in {\Bbb R}$. Then $A \in {\Bbb R}^{N \times N}$, $A$ is a normal
skew-symmetric matrix, and the eigenvalues are pure imaginary.
As an example, consider the collinear case $N = 3$. Let the particle positions be $x_1 < x_2 < x_3$, with corresponding
strengths $\Gamma_1$, $\Gamma_2$, $\Gamma_3$. The $A$ matrix is then given by
\begin{equation}
A = \left[
\begin{array}{ccc}
  0 & \frac{1}{x_1 - x_2} &  \frac{1}{x_1 - x_3}\\
 \frac{1}{x_2 - x_1} &  0 & \frac{1}{x_2 - x_3}\\
  \frac{1}{x_3 - x_1} &   \frac{1}{x_3 - x_2} &   0
\end{array}
\right].
\end{equation}
Since $N$ is odd, we have $\det(A) = 0$.
The other two eigenvalues are given by:
\begin{eqnarray}
\lambda_{123} = \pm i\sqrt{\frac{1}{(x_2 - x_1 )^2} + \frac{1}{(x_3 - x_2 )^2} + \frac{1}{(x_3 - x_1 )^2}},
\end{eqnarray}
which is invariant under cyclic permutations of the indices
($\lambda_{123} = \lambda_{231} = \lambda_{312}$).
We can scale the length of the
configuration so that the distance between $x_1$ and $x_3$ is one, hence without loss
of generality, let $x_1 = 0, x_2 = x, x_3 = 1$.
The other two eigenvalues are then given by the formula:
\begin{eqnarray}
\lambda = \pm i\sqrt{\frac{(1-x+x^2 )^2}{x^2 (1-x)^2} }.
\end{eqnarray}
It is easy to see that the numerator has no roots in the interval $(0,1)$, hence the nullspace dimension of $A$ is one.
The nullspace vector is then given (uniquely up to multiplicative constant) by:
\begin{equation}
\vec{\Gamma} =  \left[
\begin{array}{c}
1\\
 -\left( \frac{x_3 - x_2}{x_3 - x_1}\right) \\
 \left(\frac{x_3 - x_2}{x_2 - x_1}\right)
\end{array}
\right].
\end{equation}
For the special symmetric case $x_3 - x_1 = 1$, $x_3 - x_2 = 1/2$, $x_2 - x_1 = 1/2$, we have
$\Gamma_1 = 1, \Gamma_2 = -1/2, \Gamma_3 = 1$.
We show this case in figure \ref{fig2} along with the separatrices associated with the corresponding
flowfield generated by the singularities. Since the sum of the strengths of the three vortices is $\Gamma_1 + \Gamma_2 + \Gamma_3 = 1 - 1/2 + 1 = 3/2$, the far field is that of a point vortex of
strength $\Gamma = 3/2$. Interestingly, for the collinear cases, since $A$ is real, the nullspace vector is either real, or if multiplied by $i$, is pure imaginary. Hence, each collinear configuration of point vortices obtained with a given $\vec{\Gamma} \in {\Bbb R}$ is also a collinear configuration of sources/sinks with corresponding strengths given by  $i\vec{\Gamma}$. The corresponding streamline pattern for the source/sink configuration, as shown in the dashed curves of figure \ref{fig2},
is the orthogonal complement of the curves corresponding to the point vortex case.

 \begin{figure}[h]
\begin{center}
\includegraphics[scale=0.6,angle=0]{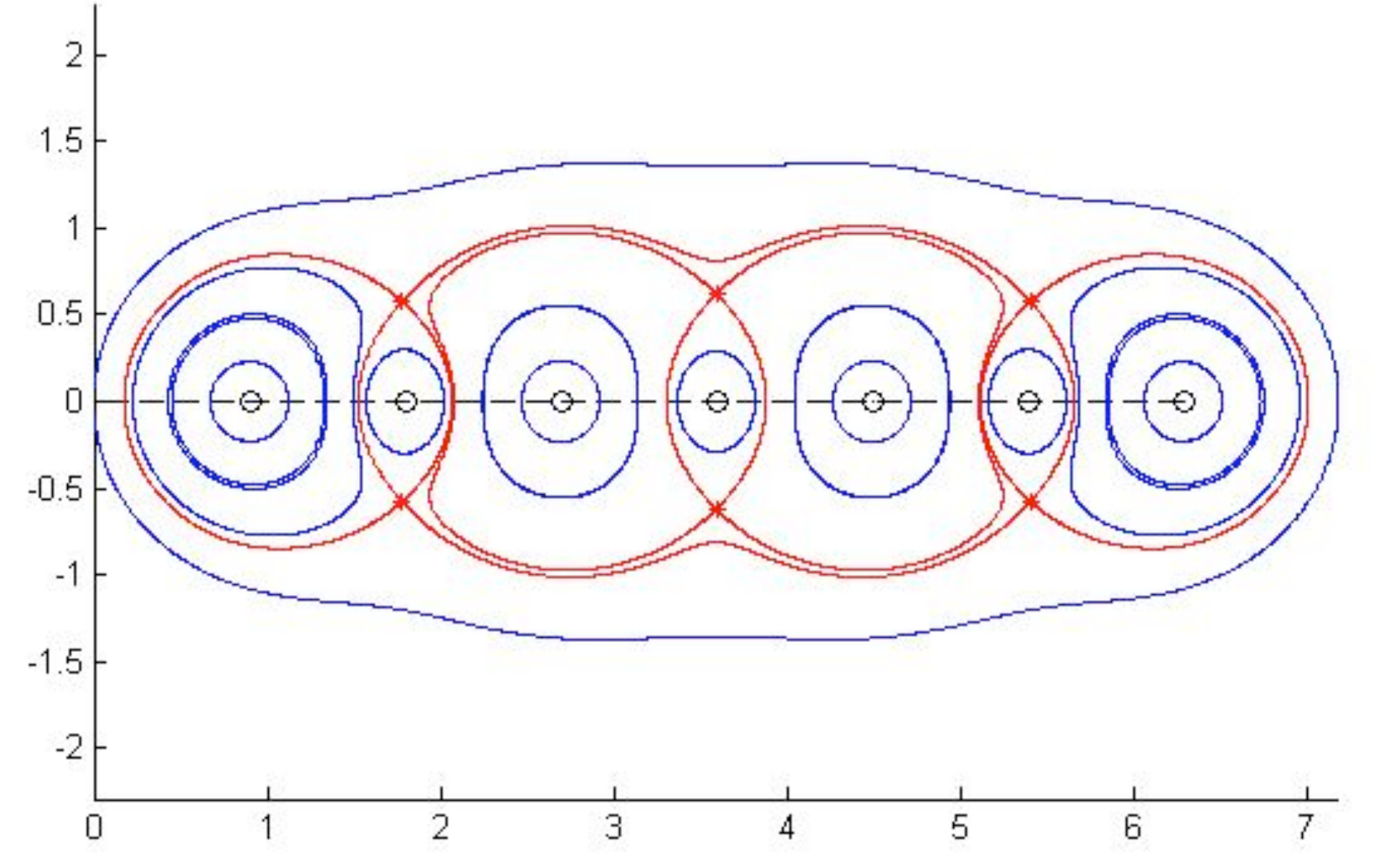}
\end{center}
\caption{\footnotesize  $N = 7$ evenly distributed point vortices on a line. The far field is that of a point vortex at the
center-of-vorticity of the system. Because of the symmetry of the spacing, the vortex strengths are symmetric about the central point $x_4$ which also corresponds to the center-of-vorticity.}
\label{fig3}
\end{figure}

 \begin{figure}[h]
\includegraphics[scale=0.6,angle=0]{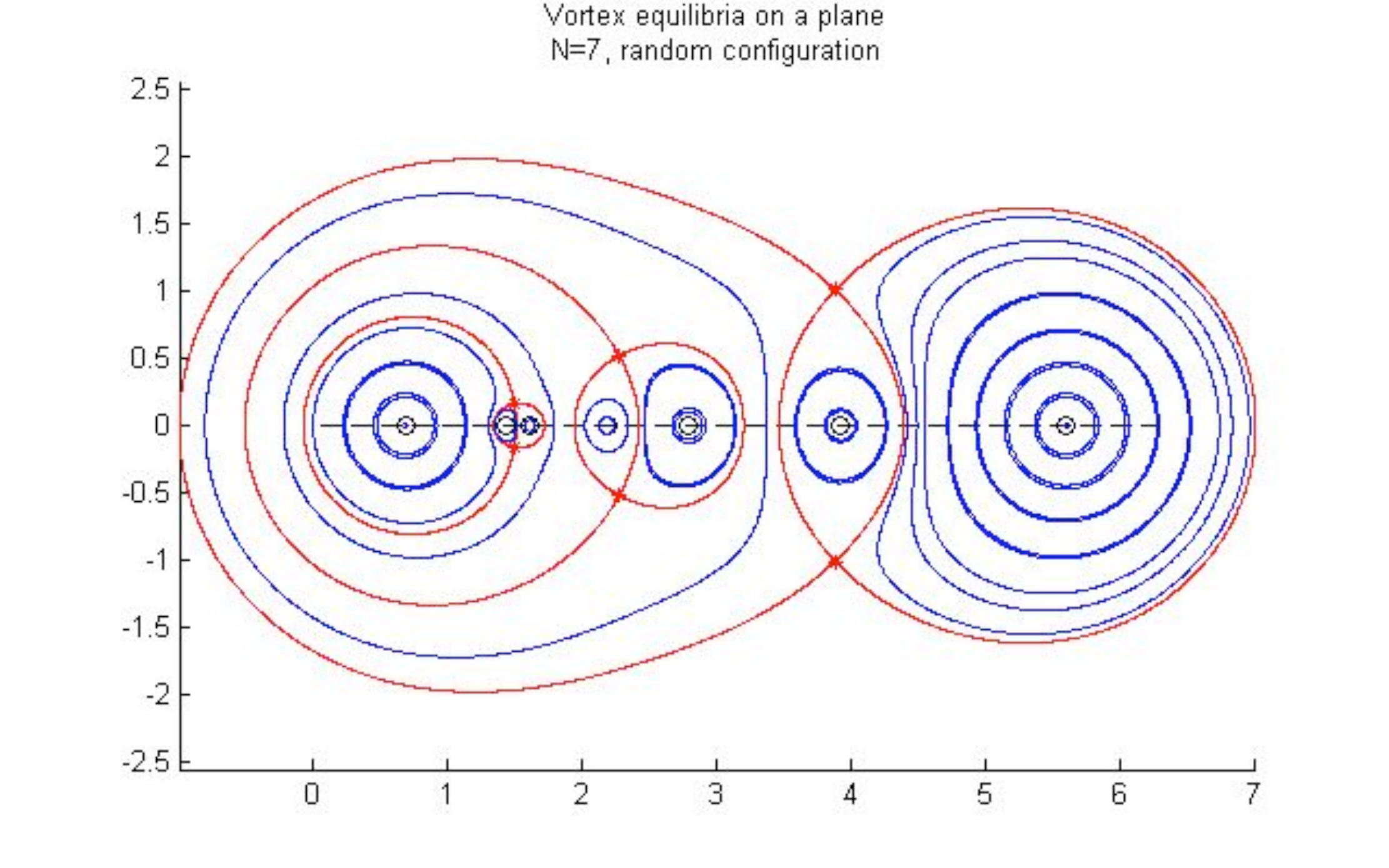}
\caption{\footnotesize
$N = 7$ randomly distributed point vortices  on a line. The far field is that of a point vortex at the
center-of-vorticity of the system.}
\label{fig4}
\end{figure}

For $N$ even, we cannot say a priori whether or not $\det (A) = 0$ as the case for $N = 2$ shows. For this, the $A$ matrix is
\begin{equation}
A = \left[
\begin{array}{cc}
  0 & \frac{1}{x_1 - x_2}\\
 \frac{1}{x_2 - x_1} &  0
\end{array}
\right] =
\left[
\begin{array}{cc}
  0 & \frac{1}{d}\\
- \frac{1}{d} &  0
\end{array}
\right].
\end{equation}
The eigenvalues are $\lambda = \pm i/d$, hence there is no equilibrium (except in the limit $d \rightarrow \infty$).

We show in figures \ref{fig3} and \ref{fig4} two representative examples of collinear fixed point vortex equilibria for $N = 7$, along with their corresponding global streamline patterns.
 In figure
\ref{fig3} we deposit seven evenly spaced points on a line and solve for the nullspace vector to obtain
the singularity strengths (ordered from left to right)
\begin{eqnarray}
\vec{\Gamma} &=& (1.0000, -0.5536, 0.9212, -0.5797, 0.9212, -0.5536, 1.0000),\\
\sum_\alpha \Gamma_\alpha &=& 2.1555.
\end{eqnarray}
Because of the even spacing, the strengths are symmetric about the central point $x_4$
($\Gamma_1 = \Gamma_7, \Gamma_2 = \Gamma_6 , \Gamma_3 = \Gamma_5$),
which is also the location of the center-of-vorticity $\sum_{\alpha = 1}^7 \Gamma_\alpha x_\alpha $.
Figure \ref{fig4} shows a fixed equilibrium corresponding to seven points randomly placed on a line.
The nullspace vector for this case is (ordered from left to right)
\begin{eqnarray}
\vec{\Gamma} &=& (1.0000, -0.5071, 0.5342, -0.4007, 0.2815, -0.2505, 1.0743),\\
\sum_\alpha \Gamma_\alpha &=& 1.7317.
\end{eqnarray}
In both cases,  the singularities are all point vortices (or source/sink systems) hence are examples of collinear equilibria such as those discussed in
Aref (2007a, 2007b, 2009) and Aref et. al. (2003)  where the strengths are typically chosen as equal.
The streamline pattern at infinity in both cases is that of a single point vortex of strength $\sum_{\alpha = 1}^7 \Gamma_\alpha \ne 0$ located at the center of vorticity $\sum_{\alpha = 1}^7 \Gamma_\alpha x_\alpha $.




\section{Triangular equilibria}

The case $N = 3$ is somewhat special and worth treating separately.
Given any three points $\{ z_1 , z_2 , z_3 \}$ in the complex plane, the corresponding configuration matrix $A$ is:

\begin{equation}
A = \left[
\begin{array}{ccc}
  0 & \frac{1}{z_1 - z_2} &  \frac{1}{z_1 - z_3}\\
 \frac{1}{z_2 - z_1} &  0 & \frac{1}{z_2 - z_3}\\
  \frac{1}{z_3 - z_1} &   \frac{1}{z_3 - z_2} &   0
\end{array}
\right].
\end{equation}
There  is no loss of generality in choosing two of the points along the real axis, one at the origin
of our coordinate system, the other at $x = 1$. Hence we set $z_1 = 0$, $z_2 = 1$, and we let
$z_3 \equiv z$.
Then $A$ is written much more simply:
\begin{equation}
A = \left[
\begin{array}{ccc}
  0 &-1& - \frac{1}{  z}\\
1&  0 & \frac{1}{1 - z}\\
  \frac{1}{z} &   \frac{1}{z - 1} &   0
\end{array}
\right].
\end{equation}
Since $N$ is odd, one of the eigenvalues of $A$ is zero. The other two are given by:
\begin{eqnarray}
\lambda = \pm i\sqrt{\frac{1}{z^2} + \frac{1}{(1-z)^2} + 1}
=
\pm i\sqrt{\frac{ (1-z+z^2 )^2}{z^2 (1-z)^2 }}
\end{eqnarray}
When the numerator is not zero, the nullspace dimension is one and it
is easy to see that the nullspace of $A$ is given by:
\begin{equation}
\vec{\Gamma} =  \left[
\begin{array}{c}
  \frac{1}{z-1} \\
- \frac{1}{z}\\
1
\end{array}
\right].\label{onednulls}
\end{equation}
However, the numerator is zero at the points:
\begin{eqnarray}
z = \exp (\frac{\pi i}{3}), \exp (\frac{5\pi i}{3}),
\end{eqnarray}
at which $\Re z = \frac{1}{2}$, $\Im z = \pm \frac{\sqrt{3}}{2}$.  This forms an equilateral triangle
in which case the nullspace dimension is three. We have thus proven the following:

\medskip

\noindent
{\bf Theorem 5.} {\it  For three point vortices, or for three sources/sinks, the only fixed equilibria are
collinear. In this case, the nullspace dimension of $A$ is one and is given by (\ref{onednulls}).
For the equilateral triangle configuration, the nullspace dimension is three.}
\medskip

\begin{figure}[h]
\begin{center}
\includegraphics[scale=0.6,angle=0]{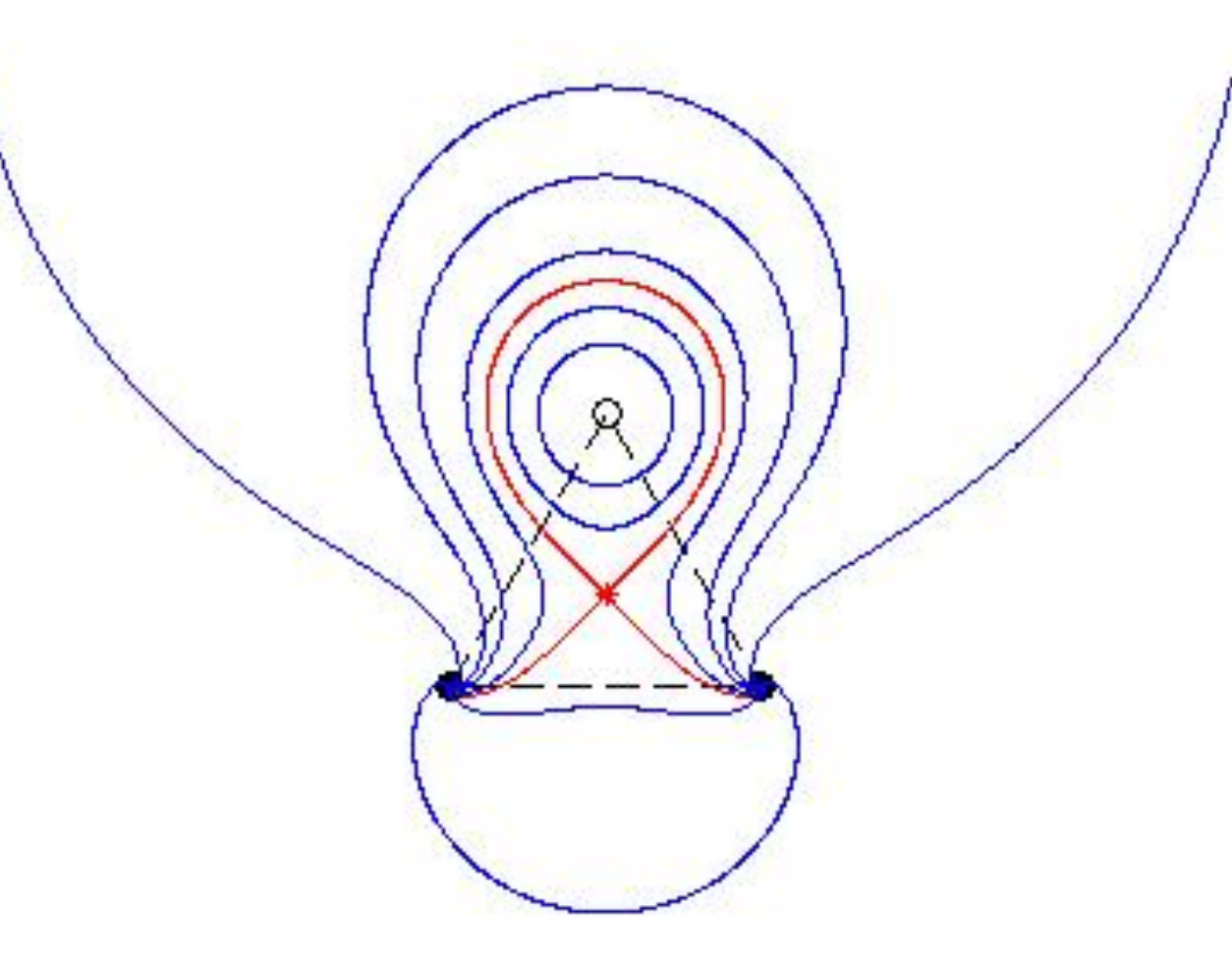}
\caption{\footnotesize $N = 3$ equilateral triangle configuration with corresponding streamline pattern.
The strengths are given by $\Gamma_1 = 1.0000$, $\Gamma_2 = -0.5000 + 0.8660i$, $\Gamma_3 = -0.5000 + 0.8660i$.
}
\label{fig5}
\end{center}
\end{figure}

\begin{figure}[h]
\begin{center}
\includegraphics[scale=0.6,angle=0]{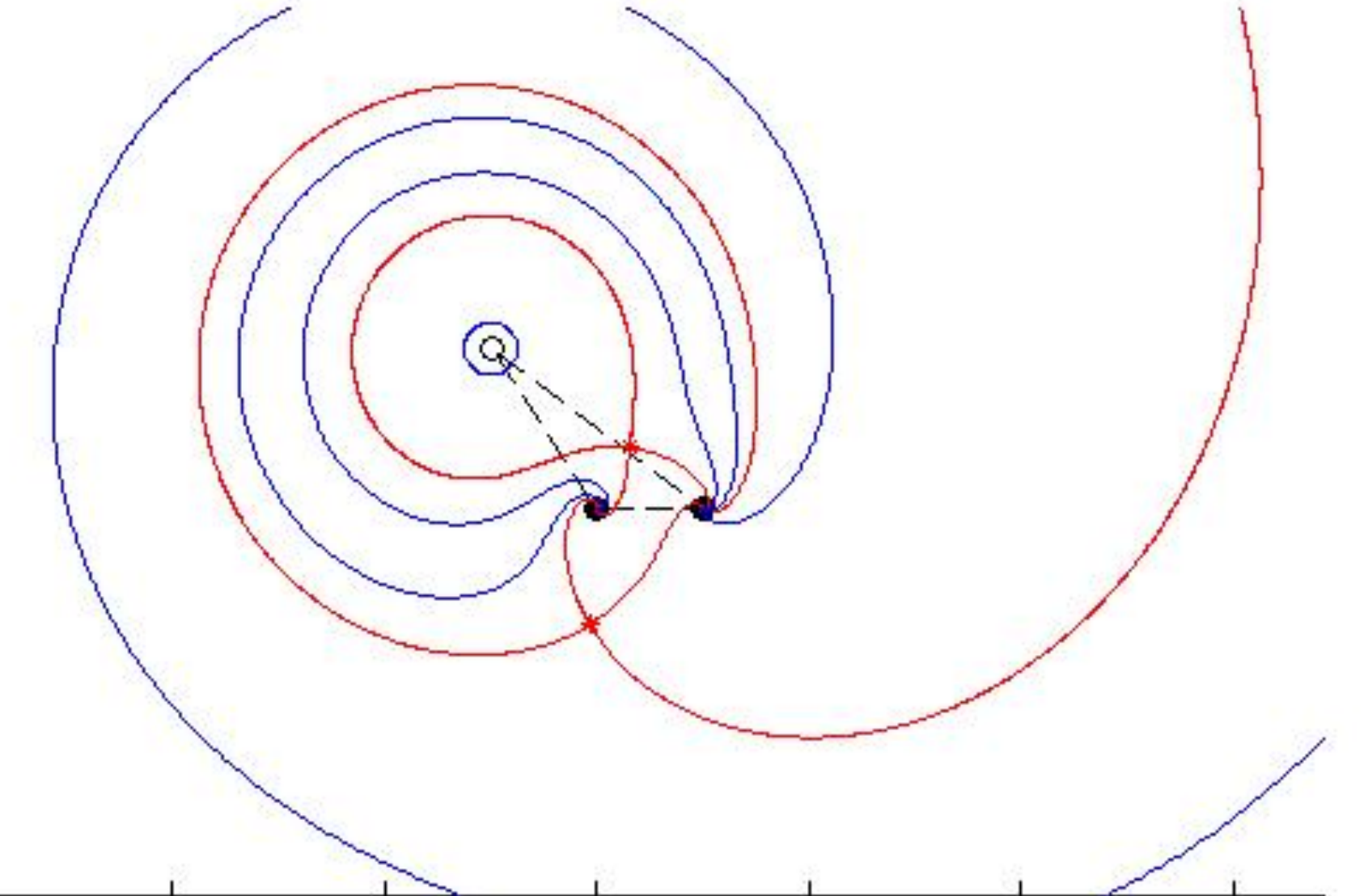}
\caption{\footnotesize $N = 3$ non-equilateral triangular state with corresponding streamline pattern.
The strengths are given by $\Gamma_1 = 1.0000$, $\Gamma_2 = 0.3077 + 0.4615i$, $\Gamma_3 = -0.3200 - 0.2400i$.
}
\label{fig6}
\end{center}
\end{figure}

\begin{figure}[h]
\begin{center}
\includegraphics[scale=0.8,angle=0]{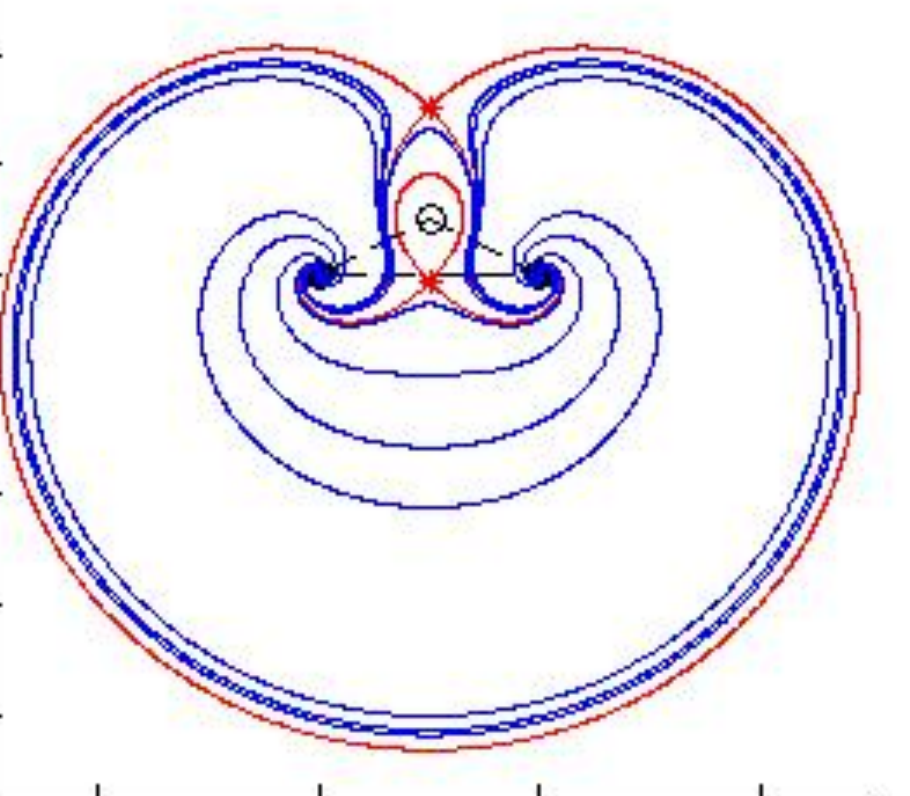}
\caption{\footnotesize  $N = 3$ isosceles triangle state.
The strengths are given by $\Gamma_1 = 1.0000$, $\Gamma_2 = -1.6000 + 0.8000i$, $\Gamma_3 = -1.6000 - 0.8000i$.
}
\label{fig7}
\end{center}
\end{figure}

\begin{figure}[h]
\begin{center}
\includegraphics[scale=0.6,angle=0]{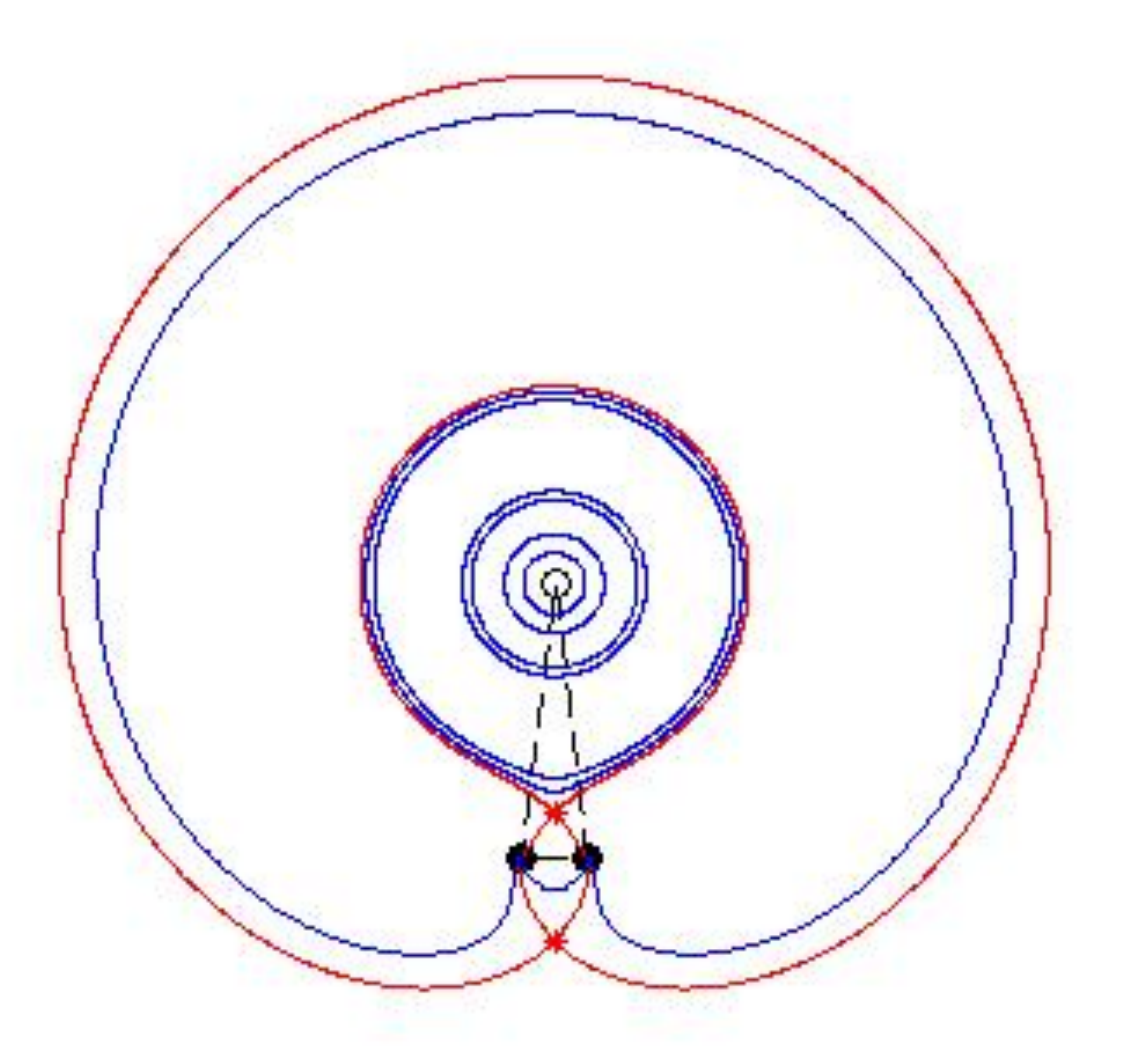}
\caption{\footnotesize $N = 3$ isosceles triangle state.
The strengths are given by $\Gamma_1 = 1.0000$, $\Gamma_2 = -0.0308 + 0.2462i$, $\Gamma_3 = -0.0308 - 0.2462i$.
}
\label{fig8}
\end{center}
\end{figure}

We show a fixed equilibrium  equilateral triangle state
 in figure \ref{fig5} along with the corresponding streamline pattern.
 Figures \ref{fig6}, \ref{fig7}, \ref{fig8} show examples of $N = 3$ triangular states that are not equilateral.

\section{Equilibria   along prescribed curves}

We now ask a more general and interesting question.  Given any curve in the complex plane, if we distribute points
$\{z_{\alpha}\}$, ($\alpha = 1, ..., N$) along the curve, is it possible to find a
 strength vector $\vec{\Gamma}$ so that the configuration is fixed? The answer is yes, if $N$ is odd,
 and sometimes, if $N$ is even.

\begin{figure}[h]
\begin{center}
\includegraphics[scale=0.75,angle=0]{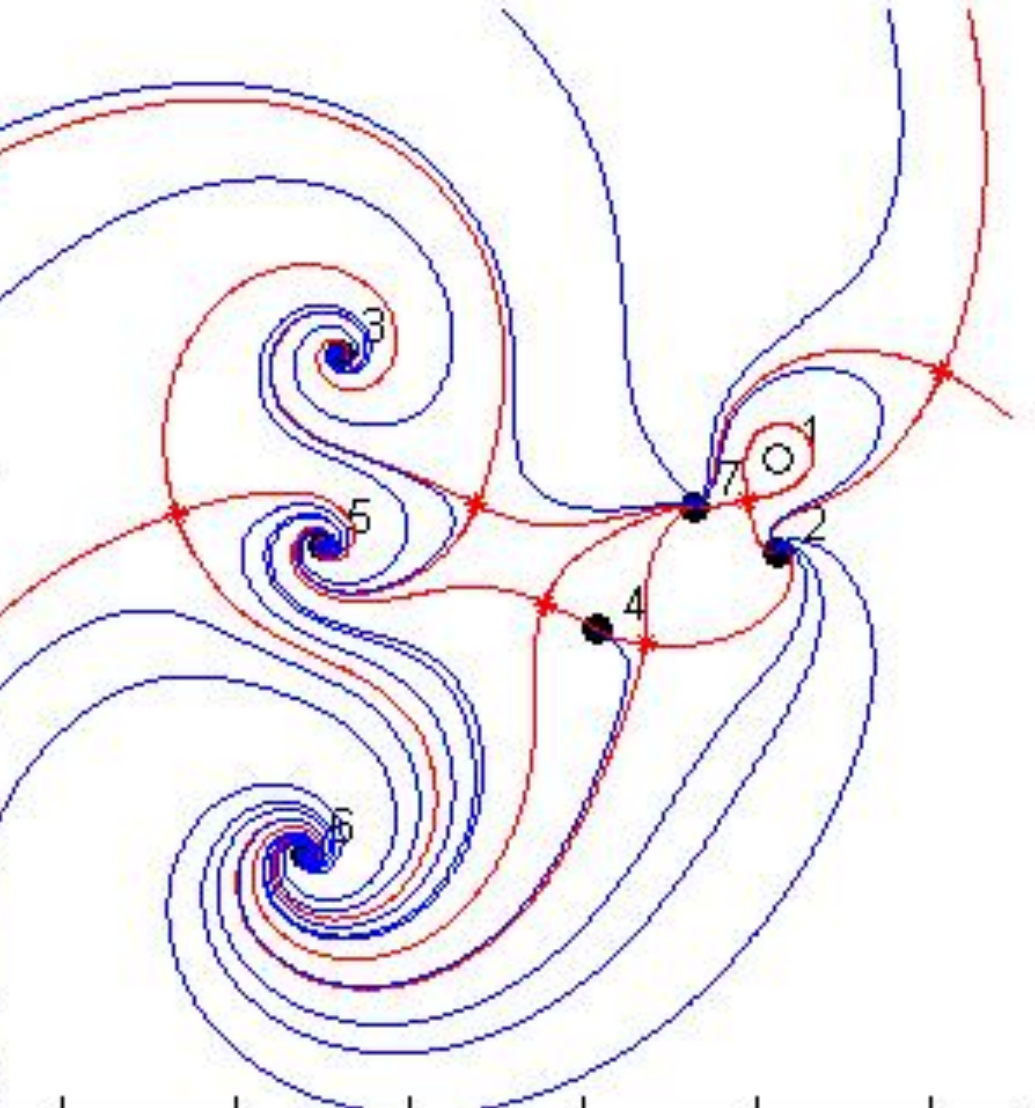}
\caption{\footnotesize Fixed equilibrium for seven points placed at random locations in the plane.
The far field is a spiral-sink (figure 1(e)) with
since $\sum \Gamma_\alpha = -2.4508 - 0.8449 i$.}
\label{fig9}
\end{center}
\end{figure}

\begin{figure}[ht]
\includegraphics[scale=0.9,angle=0]{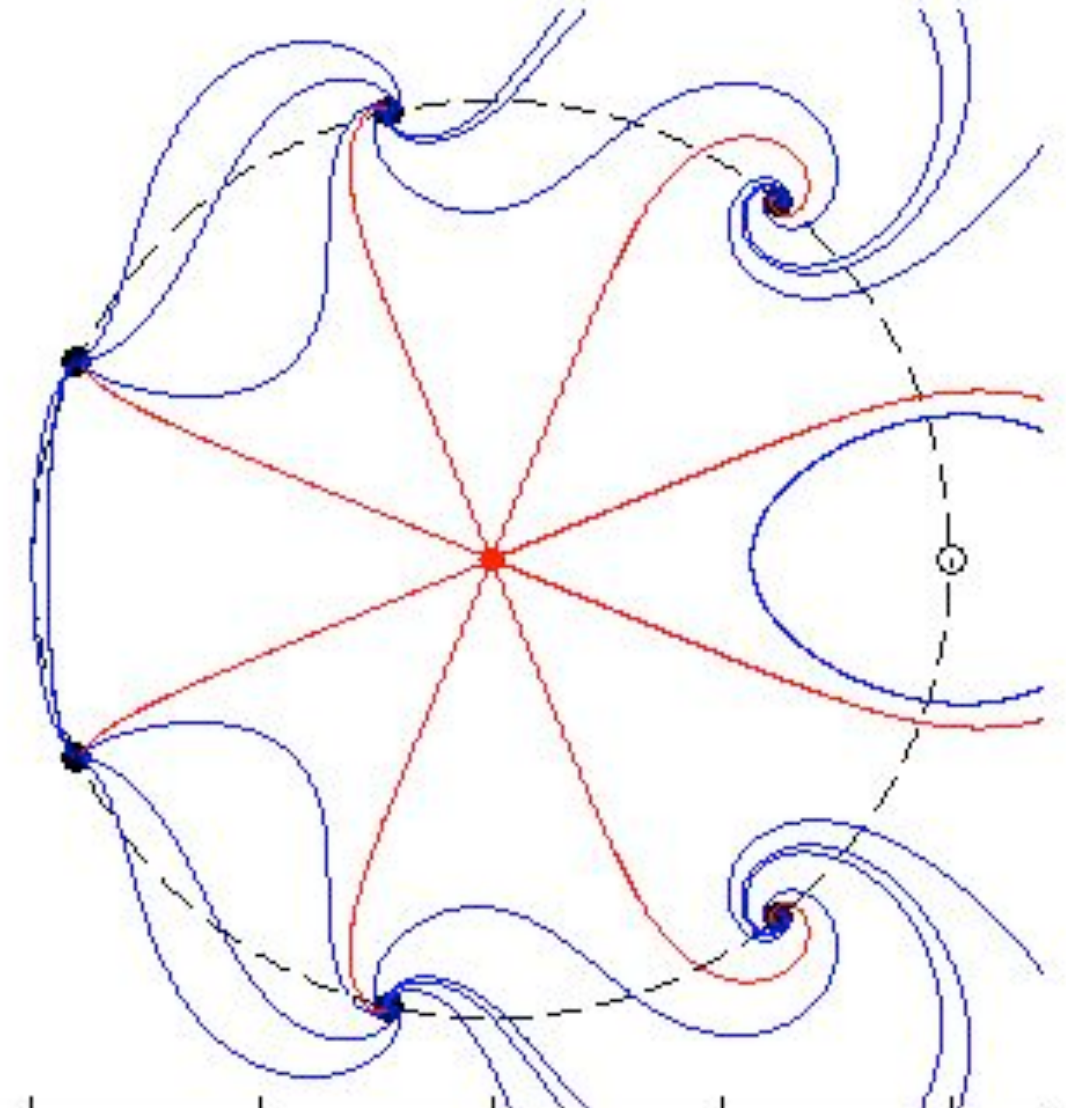}
\caption{\footnotesize  $N = 7$ evenly distributed points on a circle (dashed curve) in equilibrium. Because of the symmetry of the configuration, $\sum \Gamma_\alpha = 0$, hence the far-field vanishes.}
\label{fig10}
\end{figure}

\begin{figure}[ht]
\includegraphics[scale=0.9,angle=0]{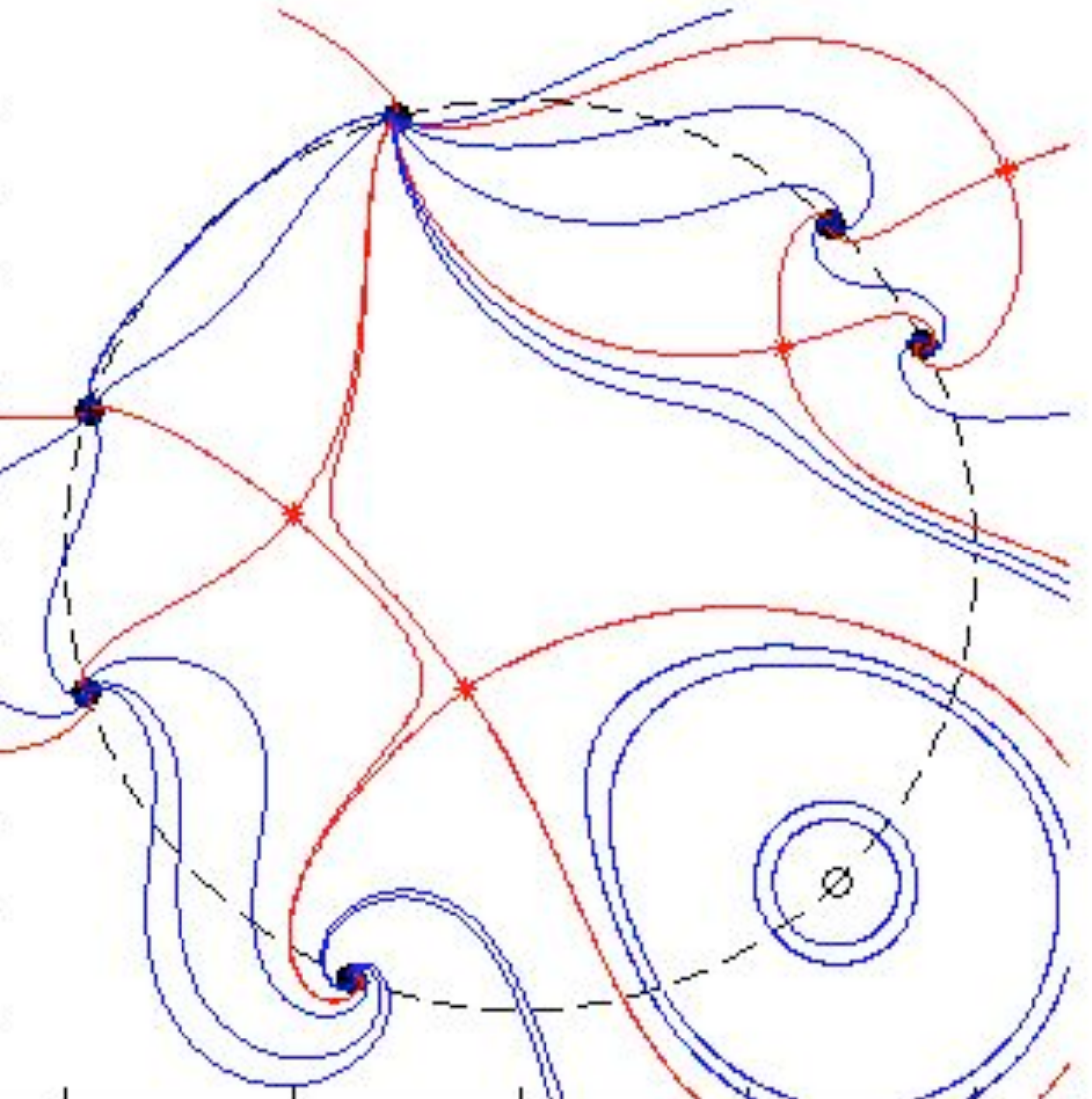}
\caption{\footnotesize
$N = 7$ randomly distributed particles on a circle (dashed curve) in equilibrium along with the corresponding streamline pattern. The far field streamline pattern is that of a spiral-sink
(figure 1(g)) since $\sum \Gamma_\alpha = 0.2649 - 0.5222 i$.}
\label{fig11}
\end{figure}

\begin{figure}[ht]
\includegraphics[scale=0.9,angle=0]{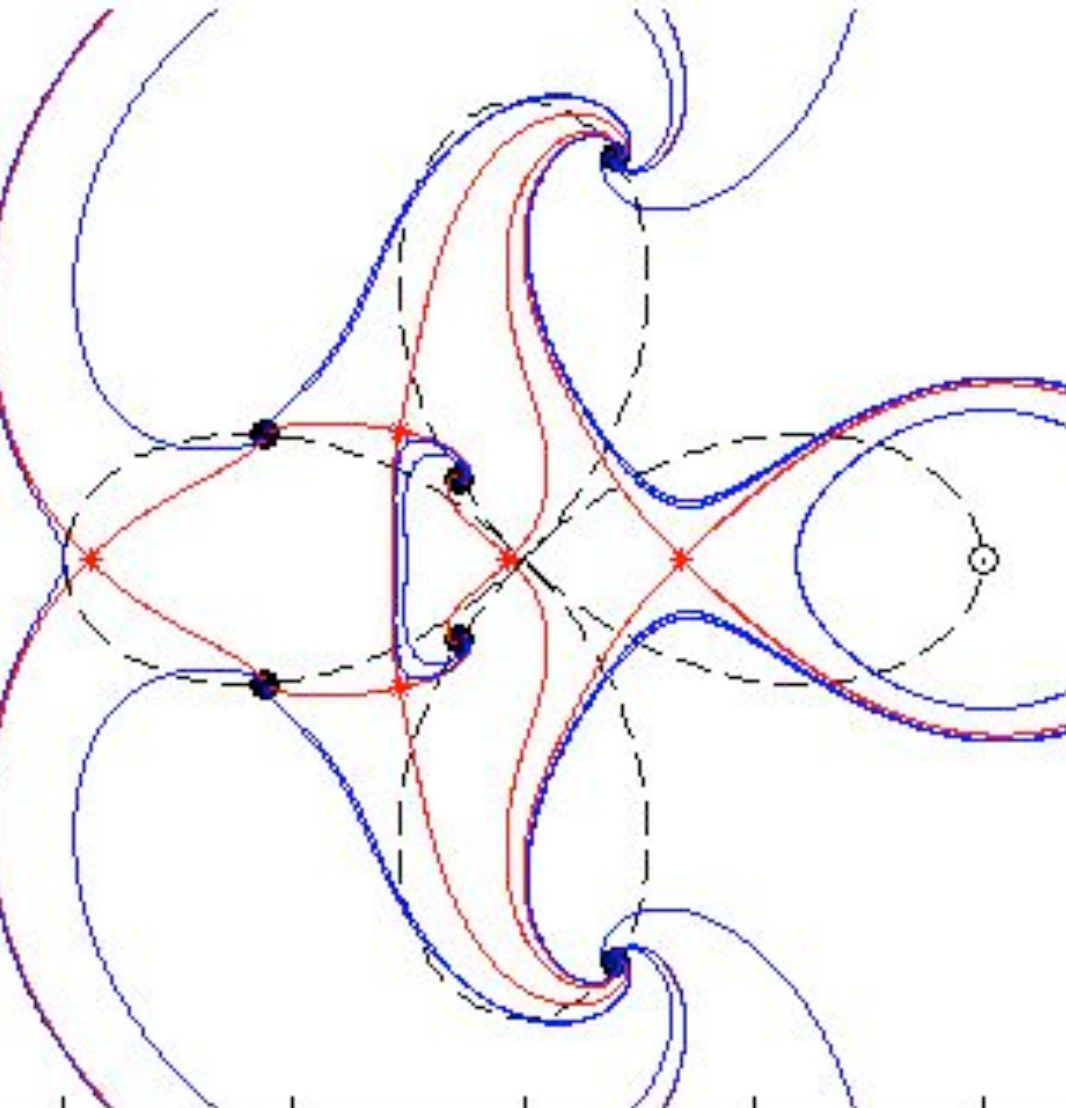}
\caption{\footnotesize
 $N = 7$ evenly   distributed particles in equilibrium  on the curve $r(\theta ) = \cos (2\theta)$
(dashed curve) along with the corresponding streamline pattern.
The far field corresponds to a point vortex since $\sum \Gamma_\alpha = -0.8892$.}
\label{fig12}
\end{figure}

\begin{figure}[ht]
\includegraphics[scale=0.9,angle=0]{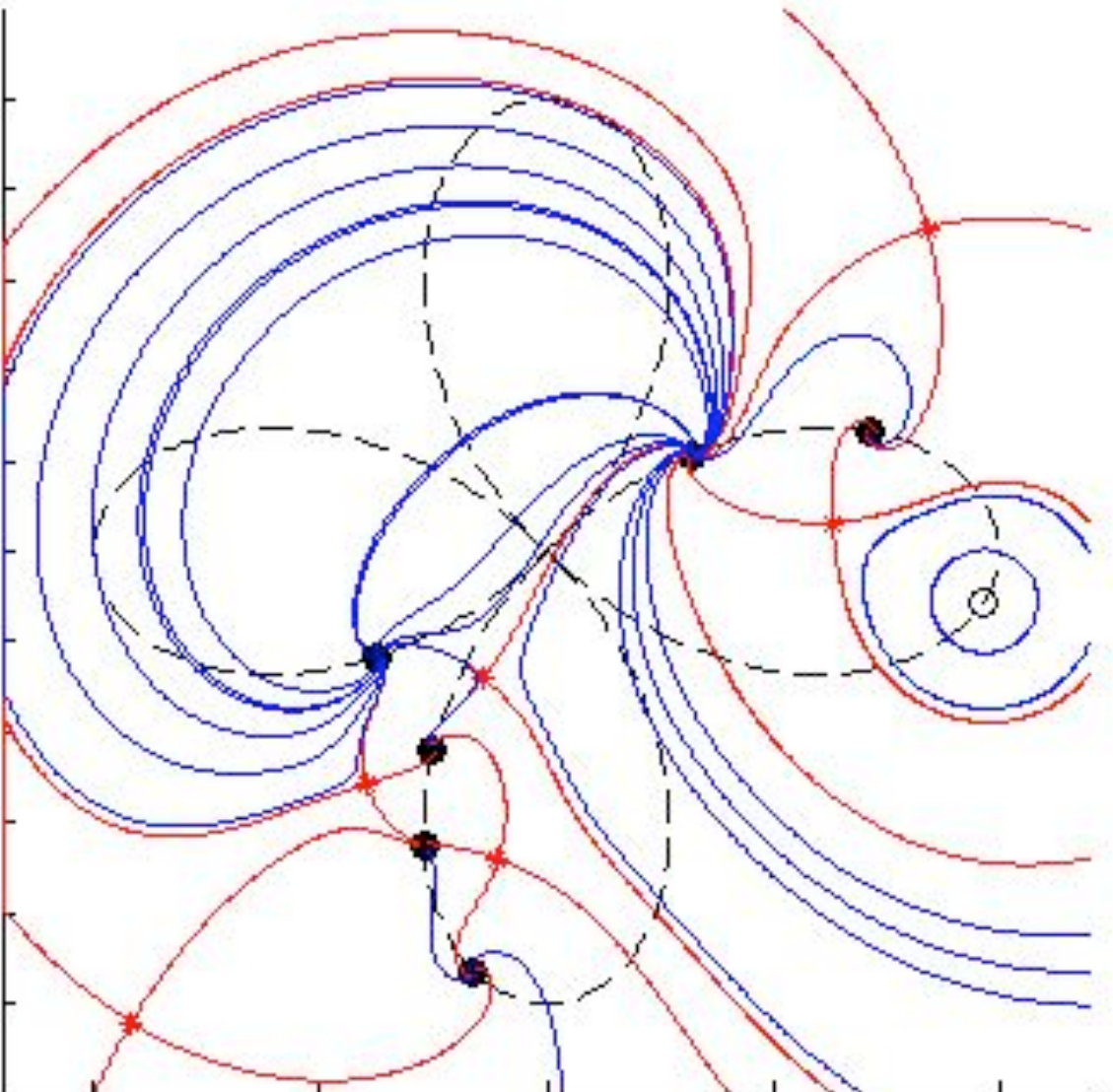}
\caption{\footnotesize
$N = 7$ randomly  distributed particles in equilibrium  on the curve $r(\theta ) = \cos (2 \theta)$
(dashed curve).
The far field corresponds to a source-spiral (figure 1(f)) since $\sum \Gamma_\alpha = 0.7244 + 0.3589i$.}
\label{fig13}
\end{figure}

\begin{figure}[ht]
\includegraphics[scale=0.9,angle=0]{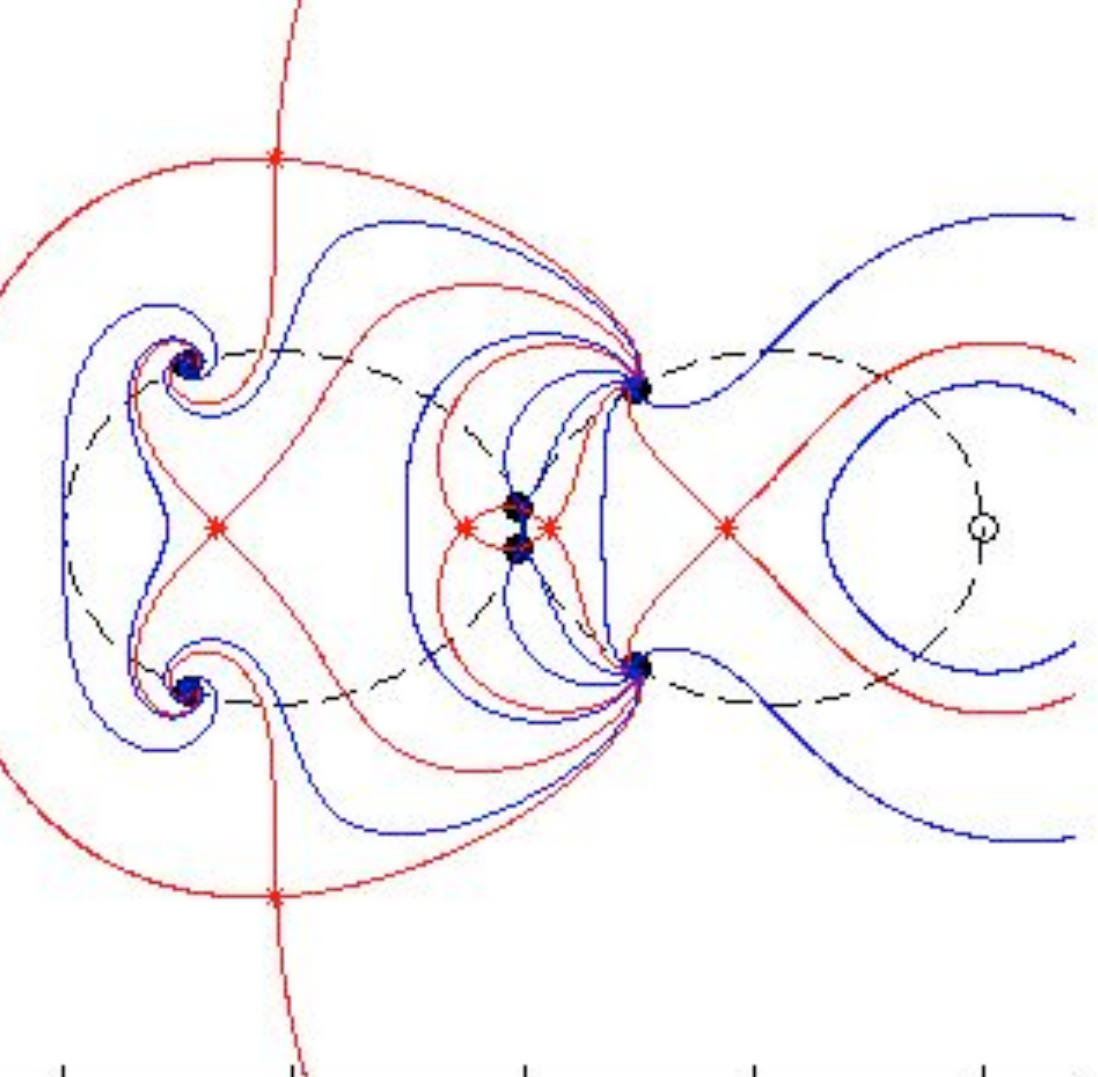}
\caption{\footnotesize $N = 7$ evenly   distributed particles in equilibrium  on the curve $r(\theta ) = \cos^2 (\theta)$ (dashed curve).
The far field corresponds to a point vortex since $\sum \Gamma_\alpha = 0.7136$.}
\label{fig14}
\end{figure}

\begin{figure}[ht]
\includegraphics[scale=0.9,angle=0]{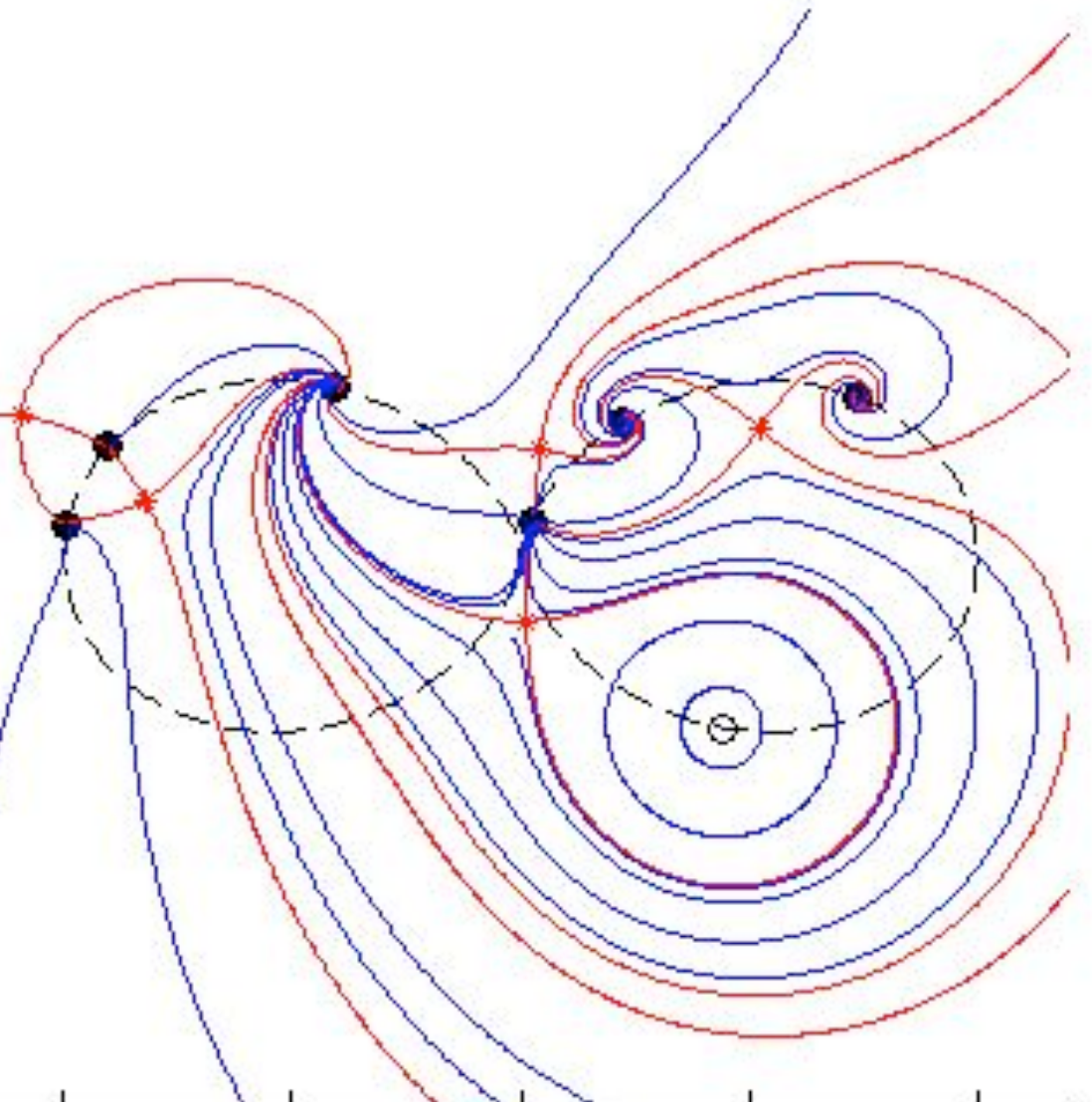}
\caption{\footnotesize
$N = 7$ randomly  distributed particles in equilibrium  on the curve $r(\theta ) = \cos^2 (\theta)$
(dashed curve).
The far field corresponds to a source-spiral (figure 1(f)) since $\sum \Gamma_\alpha = 0.9685 +1.0460i$.}
\label{fig15}
\end{figure}

Figures \ref{fig9} - \ref{fig15} show a collection of fixed equilibria along curves that we prescribe. First, figure \ref{fig9}
shows $7$ points places randomly in the plane, with the singularity strengths obtained from the nullspace of $A$ so that the
system is an equilibrium. The strengths are given by:
$\vec{\Gamma} = ( 1.0000,
 -0.7958 + 1.0089i,
 -1.3563 - 0.4012i,
  0.0297 + 0.1594i,
  0.9155 + 0.3458i,
 -2.0504 - 0.8776i,
 -0.1935 - 1.0802i)^T$
 with the sum given by $ -2.4508 - 0.8449 i$. Thus, the far field is that of a spiral-sink configuration.
 Figure \ref{fig10} shows the case of $N = 7$ points distributed evenly around a circle. The nullspace
 vector is given by $\vec{\Gamma} = (1.0000, -0.9010 + 0.4339i, 0.6235 - 0.7818i, -0.2225 + 0.9749i, -0.2225 - 0.9749i,
  0.6235 + 0.7818i, -0.9010 - 0.4339i)^T$. For this very symmetric case, the sum of the strengths is zero, hence in a sense, the
  far field vanishes.
  Figure \ref{fig11}  shows the case of $N = 7$ points placed at random positions on a circle. Here, the nullspace vector is
  given by
  $\vec{\Gamma} = ( 1.0000,
 -0.6342 + 0.4086i,
  0.3699 - 0.5929i,
 -0.1501 + 0.6135i,
 -0.2483 - 0.9884i,
  0.2901 + 0.3056i,
 -0.3595 - 0.2686i)^T$ The random placement of points breaks the symmetry of the previous case and the sum of strengths is given by $0.2649 - 0.5222 i$ which corresponds to a spiral-sink. In figures \ref{fig12} and \ref{fig13} we show equilibrium
 distribution of points along a curve we call a `flower-petal', given by the formula $r(\theta ) = \cos (2\theta)$, $0 \le \theta \le 2 \pi$.
 In figure \ref{fig12} we distribute them evenly on the curve, while in figure \ref{fig13} we distribute them randomly.
The particle strengths from the configuration in figure \ref{fig12} are $\vec{\Gamma} = (1.0000,
 0.1824 + 0.1498i,
 -0.9892 - 0.9103i,
 -0.1378 - 0.5333i,
 -0.1378 + 0.5333i,
 -0.9892 + 0.9103i,
  0.1824 - 0.1498i )^T $ with sum equaling $ -0.8892$ corresponding to a far field point vortex. Figure \ref{fig13} shows particles
  distributed randomly on the same flower-petal curve. Here, the particle strengths are $\vec{\Gamma} = (1.0000,
0.2094 - 0.4071i,
 -0.3009 + 0.3003i,
0.0404 - 0.2864i,
 -0.1779 + 0.2773i,
 0.4236 + 0.8052i,
 -0.4702 - 0.3304i)^T$, with sum given by $ .7244 + .3589i$. Hence the far field corresponds to a source-spiral.

The last two configurations, shown in figures \ref{fig14} and \ref{fig15} are equilibria distributed along figure eight curves, given by the formulas
$r(\theta ) = \cos^2 (\theta)$, $0 \le \theta \le 2 \pi$. In figure \ref{fig14}  we distribute the points evenly around the curve, which
gives rise to strengths $\vec{\Gamma} = (1.0000,
 -0.2734 + 0.5350i,
  0.0239 - 0.2080i,
  0.1063 - 0.0517i,
  0.1063 + 0.0517i,
  0.0239 + 0.2080i,
 -0.2734 - 0.5350i )^T $, whose sum is $.7136$, thus a far field point vortex. In contrast, when the points are distributed
 randomly around the same curve, as in figure \ref{fig15}, the strengths are given by
 $\vec{\Gamma} = (1.0000, -0.1054 + 0.5724i,  -0.0174 - 0.4587i,   0.9208 + 1.2450i,  -0.0460 - 0.4577i,  -0.5292 + 0.2371i,  -0.2543 - 0.0921i)^T $, with sum equaling $ .9685+1.0460i$, hence a far field source-spiral.


\section{Classification of equilibria in terms of the singular spectrum}

 \begin{figure}[ht]
\includegraphics[scale=0.7,angle=0]{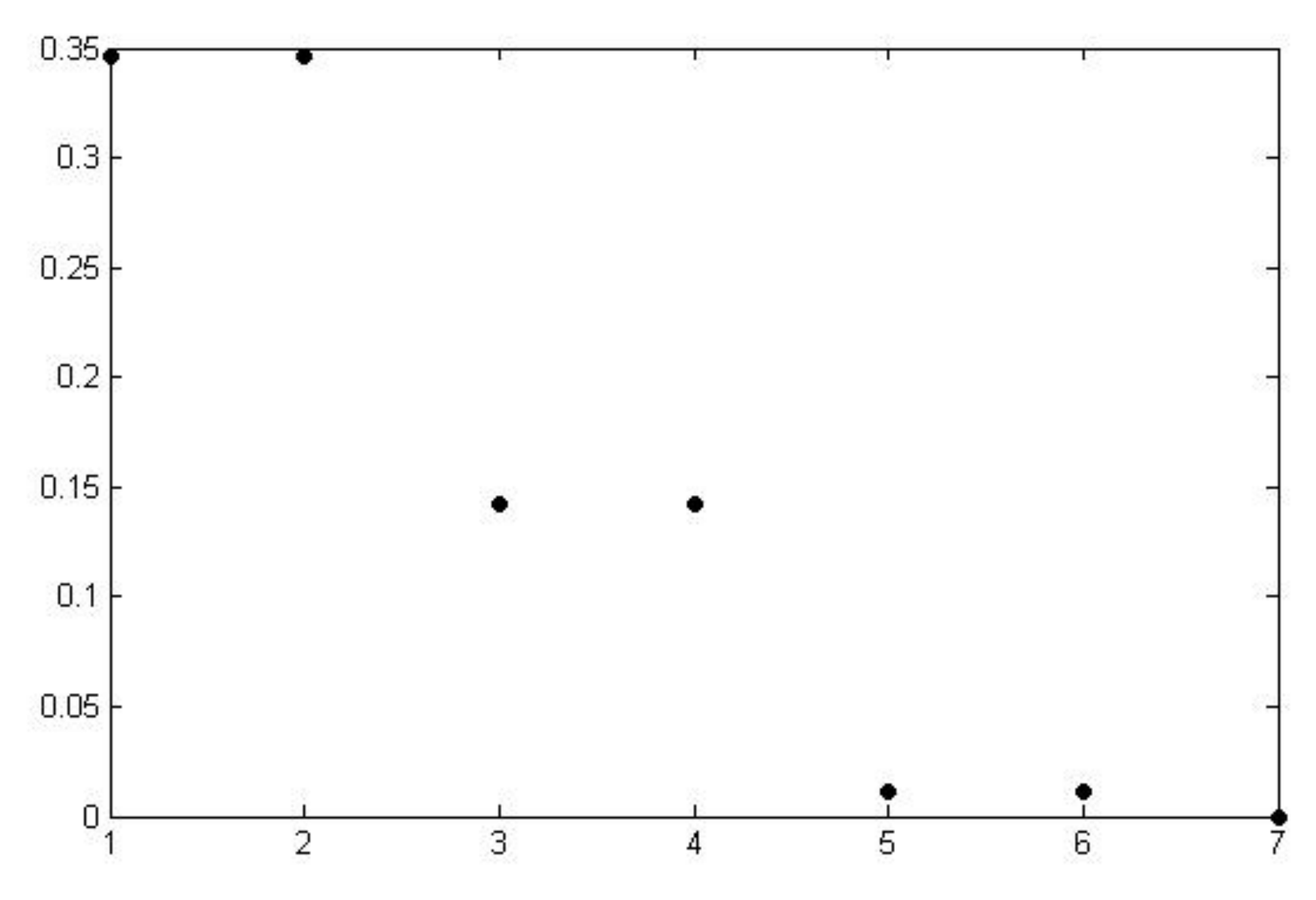}
\caption{\footnotesize Singular spectrum (normalized) for $N = 7$ particles placed randomly along a figure eight planar curve
(i.e. the equilibrium configuration shown in figure \ref{fig15}).
Singular values are grouped in pairs (except for the zero value) due to the skew-symmetry of the configuration matrix.}
\label{fig16}
\end{figure}

\begin{table}
\begin{center}
 \begin{minipage}{8.0cm}
\begin{tabular}{||l|c|c|c||}
\hline
      Configuration & $\sigma$ (unormalized)  &  $\sigma$ (normalized)  & Shannon entropy  \\\hline
      & 1.0000  & 0.5000 & 0.6931\\
    Equilateral  & 1.0000 & 0.5000 &  \\
      & 0.00 & 0.00 &
      \\\hline
      &1.0598 & 0.5000 & 0.6931  \\
    Isosceles (acute)& 1.0598 & 0.5000 &  \\
      & 0.00 & 0.00 &
      \\\hline
       &2.7203 & 0.5000 & 0.6931  \\
     Isosceles (obtuse)& 2.7203 & 0.5000 &  \\
      & 0.00 & 0.00 &
      \\\hline
       &1.2115 & 0.5000 & 0.6931  \\
   Arbitrary triangle & 1.2115 & 0.5000 &  \\
      & 0.00 & 0.00  &
       \\\hline
 \end{tabular}
  \end{minipage}
\end{center}
\caption{Singular spectrum  of triangular states ($N = 3$)
\label{table1}}
\end{table}

\begin{table}
\begin{center}
 \begin{minipage}{8.0cm}
\begin{tabular}{||l|c|c|c||}
\hline
      Configuration & $\sigma$ (unormalized)  &  $\sigma$ (normalized)  & Shannon entropy  \\\hline
      & 4.5000  & 0.5000 & 0.6931\\
     $N=3$  & 4.5000 & 0.5000 &  \\
      & 0.00 & 0.00 &
      \\\hline
      & 2.5249  & 0.3214 & 1.5237  \\
    $N = 7$ (even) & 2.5249 & 0.3214 &  \\
      & 1.6831 & 0.1428 & \\
       & 1.6831 & 0.1428 &\\
     & 0.8420 & 0.0357 & \\
      & 0.00  & 0.00 &
      \\\hline
       & 6.3408  & 0.4457 & 1.0723  \\
    $N = 7$ (random) & 6.3408 & 0.4457 &  \\
      & 2.0969 & 0.0487 & \\
       & 2.0969 & 0.0487 &\\
     & 0.7062 & 0.0055 & \\
      & 0.7062  & 0.0055 & \\
       & 0.0000 & 0.0000 &
      \\\hline
        \end{tabular}
  \end{minipage}
\end{center}
\caption{Singular spectrum of collinear states ($N = 3, 7$)
\label{table2}}
\end{table}

\begin{table}
\begin{center}
 \begin{minipage}{8.0cm}
\begin{tabular}{||l|c|c|c||}
\hline
      Configuration & $\sigma$ (unormalized)  &  $\sigma$ (normalized)  & Shannon entropy  \\\hline
      & 3.0000  & 0.3214 & 1.5236  \\
    $N = 7$ (even) & 3.0000 & 0.3214 &  \\
      & 2.0000 & 0.1429 & \\
       & 2.0000 & 0.1429 &\\
     & 1.0000 & 0.0357 & \\
      & 1.0000  & 0.0357 &  \\
      & 0.0000 & 0.0000 &
      \\\hline
       & 3.7954  & 0.3363 & 1.4700 \\
    $N = 7$ (random) & 3.7954 & 0.3363 &  \\
      & 2.4250 & 0.1373 & \\
       & 2.4250 & 0.1373 &\\
     & 1.0631 & 0.0264 & \\
      & 1.0631  & 0.0264 & \\
       & 0.0000 & 0.0000 &
      \\\hline
        \end{tabular}
  \end{minipage}
\end{center}
\caption{Singular spectrum of circular states ($N = 7$)
\label{table3}}
\end{table}

\begin{table}
\begin{center}
 \begin{minipage}{8.0cm}
\begin{tabular}{||l|c|c|c||}
\hline
      Configuration & $\sigma$ (unormalized)  &  $\sigma$ (normalized)  & Shannon entropy  \\\hline
      & 11.9630  & 0.4664 & 0.9651  \\
    $N = 7$ (even) & 11.9630 & 0.4664 &  \\
      & 3.0001 & 0.0293 & \\
       & 3.0001 & 0.0293 &\\
     & 1.1454 & 0.0043 & \\
      & 1.1454  & 0.0043 &  \\
      & 0.0000 & 0.0000 &
      \\\hline
       & 6.9337  & 0.3465 & 1.3929 \\
    $N = 7$ (random) & 6.9337 & 0.3465 &  \\
      & 4.4357 & 0.1418 & \\
       & 4.4357 & 0.1418 &\\
     & 1.2769 & 0.0117 & \\
      & 1.2769  & 0.0117 & \\
       & 0.0000 & 0.0000 &
      \\\hline
        \end{tabular}
  \end{minipage}
\end{center}
\caption{Singular spectrum of figure eight states ($N = 7$)
\label{table4}}
\end{table}

\begin{table}
\begin{center}
 \begin{minipage}{8.0cm}
\begin{tabular}{||l|c|c|c||}
\hline
      Configuration & $\sigma$ (unormalized)  &  $\sigma$ (normalized)  & Shannon entropy  \\\hline
      & 5.9438  & 0.4447 & 1.1034  \\
    $N = 7$ (even) & 5.9438 & 0.4447 &  \\
      & 1.8115 & 0.0413 & \\
       & 1.8115 & 0.0413 &\\
     & 1.0538 & 0.0140 & \\
      & 1.0538  & 0.0140 &  \\
      & 0.0000 & 0.0000 &
      \\\hline
       & 8.0780  & 0.3875 & 1.3393 \\
    $N = 7$ (random) & 8.0780 & 0.3875 &  \\
      & 3.8900 & 0.0899 & \\
       & 3.8900 & 0.0899 &\\
     & 1.9523 & 0.0226 & \\
      & 1.9523  & 0.0226 & \\
       & 0.0000 & 0.0000 &
      \\\hline
        \end{tabular}
  \end{minipage}
\end{center}
\caption{Singular spectrum of flower states ($N = 7$)
\label{table5}}
\end{table}

Here we describe how to
use the non-zero singular spectrum of $A$ to classify the equilibrium when $A$ has a kernel.
Let $\sigma^{(i)}$, $i = 1, ..., k < N$ denote the non-zero singular values of the configuration matrix $A$, arranged in
descending order $\sigma^{(1)} \ge \sigma^{(2)} \ge ... \ge \sigma^{(k)} > 0$.  First we normalize each of the singular values so that they sum to one:

\begin{eqnarray}
\hat{\sigma}^{(i)} \equiv \sigma^{(i)} /\sum_{j=1}^k \sigma^{(j)}
\end{eqnarray}
Then
\begin{eqnarray}
\sum_{i=1}^k \hat{\sigma}^{(i)} =1,
\end{eqnarray}
and the string of $k$ numbers arranged from largest to smallest: $(\hat{\sigma}^{(1)} , \hat{\sigma}^{(2)}, ..., \hat{\sigma}^{(k)})$
is the `spectral representation' of the equilibrium.
The rate at which they decay from largest to smallest
is encoded in a scalar quantity called the Shannon entropy, $S$, of the matrix (see Shannon (1948) and more recent discussions associated with vortex lattices in
Newton and Chamoun (2009)):

\begin{eqnarray}
S = - \sum_{i=1}^k \hat{\sigma}^{(i)} \log \hat{\sigma}^{(i)} .
\end{eqnarray}
With this representation, spectra that drop off rapidly from highest to lowest, are `low-entropy equilibria', whereas those
that drop off slowly (even distribution of normalized singular values) are `high-entropy equilibria'.
Note that from the representation (\ref{svddecomp}), low-entropy equilibria have configuration matrix representations that are dominated
in size by a small number of terms, whereas the configuration matrices of high-entropy equilibria equilibria have terms
that are more equal in size.
See Newton and Chamoun (2009) for more detailed discussions in the context of relative equilibrium configurations, and the original report of Shannon (1948) which has illuminating discussions of entropy, information content, and its interpretations
with respect to randomness.
As an example of the normalized spectral distribution associated with the figure eight equilibrium shown in figure \ref{fig15},
we show in figure \ref{fig16} the $7$ singular values (including the zero one). The fact that they are grouped in pairs follows
from the skew-symmetry of $A$ which implies that the eigenvalues come in pairs $\pm \lambda$. Since the singular values
are the squares of the eigenvalues, it follows that there are two of each of the non-zero ones.

Tables \ref{table1} - \ref{table5} show the complete singular spectrum for all the equilibria considered in this paper.
A common measure of `robustness' associated with the configuration matrix, hence the equilibrium, is the size of the
`spectral gap' as measured by
 the size of the smallest non-zero singular value.
  From Table \ref{table2}, the collinear state with points distributed randomly and the figure-eight state with points distributed evenly (Table \ref{table4}) are the least robust in that their smallest
non-zero singular values are closest to zero.

\section{Discussion}

In this paper we describe a new method for finding and classifying   fixed equilibrium distributions
of point singularities of source/sink, point vortex, or spiral source/sink type in the complex plane under the dynamical assumption that each point
`goes with the flow'.  This includes configurations placed at random points in the
plane, at prescribed points, or lying along prescribed curves. This last situation is reminiscent of
a  classical technique for enforcing
boundary conditions along arbitrarily shaped boundaries embedded in fluid flows. These techniques are
generically referred to as singularity distribution methods. See, for example,
Katz and Plotkin (2001) and Cortez (1996, 2000) for applications and discussions of these methods in the context of
potential flow, hence inviscid boundary conditions, and Cortez (2001)  in the context of Stokes flow, hence viscous boundary conditions.
For these problems, there is generally no associated evolution equation for the interacting singularities which discretize the boundary, as in (\ref{eqn1a}) for us.
Their positions are fixed to lie along the given boundary, and the strengths  are then judiciously  chosen to enforce the
relevant inviscid or viscous boundary conditions. As in Cortez (2001), it would  be of interest to `regularize' the point singularities (1.1) and ask if the methods in this paper can be extended to smoothed out singularities, as
would a more complete analysis of the `pseudo-spectrum' associated with the configuration matrices $A$, as discussed  Trefethen and Bau (1997).

{}
\end{document}